# The Poisson formula for groups with hyperbolic properties

By Vadim A. Kaimanovich

**Abstract**

The Poisson boundary of a group $G$ with a probability measure $\mu$ is the space of ergodic components of the time shift in the path space of the associated random walk. Via a generalization of the classical Poisson formula it gives an integral representation of bounded $\mu$-harmonic functions on $G$. In this paper we develop a new method of identifying the Poisson boundary based on entropy estimates for conditional random walks. It leads to simple purely geometric criteria of boundary maximality which bear hyperbolic nature and allow us to identify the Poisson boundary with natural topological boundaries for several classes of groups: word hyperbolic groups and discontinuous groups of isometries of Gromov hyperbolic spaces, groups with infinitely many ends, cocompact lattices in Cartan-Hadamard manifolds, discrete subgroups of semi-simple Lie groups.

## 0. Introduction

The classical Poisson integral representation formula for the harmonic functions on the hyperbolic plane $\mathbb{H}^2$ can be written as

$$(0.1) \qquad f(x) = \langle F, \nu_x \rangle$$

where $\nu_x$ are the harmonic measures on the circle at infinity $\partial \mathbb{H}^2$ associated with the points $x \in \mathbb{H}^2$. The right-hand side of (0.1) makes sense for any bounded measurable function $F \in L^\infty(\partial \mathbb{H}^2)$, and (0.1) establishes an isometry between the Banach space of bounded harmonic functions $f$ on $\mathbb{H}^2$ and the space $L^\infty(\partial \mathbb{H}^2)$. [Throughout this introduction the reader is referred to the author's survey [Ka96] for all historical and background references.]

In fact, the *Poisson formula* (0.1) can be put into the much more general context of the theory of Markov operators (recall that the classical harmonic functions are characterized by a mean value property with respect to the heat kernel). For any Markov operator $P$ on a Lebesgue measure space $(X, m)$ there exists a space $\Gamma$ (the *Poisson boundary* of $P$) endowed with a family of probability measures $\nu_x$, $x \in X$ such that the Poisson formula (0.1) establishes an isometry between the space of *P-harmonic functions* (i.e., such that $Pf = f$) from $L^\infty(X, m)$ and the space $L^\infty(\Gamma)$. The Poisson boundary is defined as the space of ergodic components of the time shift $T$ in the space of sample paths of the Markov chain on $X$ associated with the operator $P$, the measures $\nu_x$ being the images of the measures in the path space corresponding to starting the Markov chain from points $x \in X$.

We emphasize that the Poisson boundary is a purely measure theoretical object (unlike the topological *Martin boundary*). If the Martin boundary is well-defined, then, viewed as a measure space with the representing measure of the constant harmonic function, it is isomorphic to the Poisson boundary.

By definition, the Poisson boundary is the *maximal* among all the spaces $B$ such that there exists a measurable map $\Pi$ from the path space to $B$ with the property that

$$(0.2) \qquad \Pi(\boldsymbol{x}) = \Pi(T\boldsymbol{x}) \qquad \text{for a.e. sample path } \boldsymbol{x} = \{x_n\}$$

(one can say that the Poisson boundary completely describes the stochastically significant behaviour of the sample paths at infinity). An example of such a space (*a priory* not maximal!) arises in the situation when $X$ is embedded into a topological space $\overline{X}$, and a.e. sample path $\boldsymbol{x}$ converges to a limit $\Pi\boldsymbol{x} = \lim x_n \in \overline{X}$.

Often, the state space $X$ is endowed with additional (geometrical, combinatorial, algebraic, etc.) structures, and the operator $P$ complies with (or, is governed by) them. Then a natural question is to *identify* (describe) the Poisson boundary in terms of these structures. This problem usually splits into two quite different parts:

(1) to find a space $B$ (related to the structure of $X$) and a map $\Pi$ with the property (0.2). *A priori* such a space is just a quotient of the Poisson boundary,

(2) to prove that the space $B$ is in fact maximal, i.e., is isomorphic to the whole Poisson boundary.

In other words, first one has to exhibit a certain system of invariants ("patterns") of the behaviour of the Markov chain at infinity, and then show completeness of this system, i.e., that these patterns completely describe the behaviour at infinity.

In the present paper we address the problem of identification of the Poisson boundary for the *random walk* determined by a probability measure $\mu$ on a countable group $G$. The transitions of the random walk are $x \mapsto xh$, where the increment $h$ is $\mu$-distributed, so that the random walk is a Markov chain homogeneous both in time and in space. In this situation the Poisson boundary is endowed with a natural action of the group $G$, and for the purposes of the identification of the Poisson boundary one may consider its equivariant quotients only (they are called *$\mu$-boundaries*).

We develop a new method (announced in the author's notes [Ka85], [Ka94]) of proving maximality of $\mu$-boundaries based on entropies of conditional random walks. Denote by $\mathbf{P}$ the measure in the path space $G^{\mathbb{Z}_+}$ corresponding to starting the random walk from the group identity. Since a $\mu$-boundary $B$ is a quotient of the path space, the points of $B$ determine conditional measures of $\mathbf{P}$. These measures correspond to the Markov chains on $G$ (*conditional random walks*) which are still homogeneous in time but lose the spatial homogeneity. Let us say that a probability measure $\Lambda$ in the path space $G^{\mathbb{Z}_+}$ has *asymptotic entropy* $\mathbf{h}(\Lambda)$ if its one-dimensional distributions $\lambda_n$ have the following Shannon-McMillan-Breiman type equidistribution property: $-\log \lambda_n(x_n) \to \mathbf{h}(\Lambda)$ for $\Lambda$-a.e. $\boldsymbol{x} = \{x_n\} \in G^{\mathbb{Z}_+}$ and in the space $L^1(\Lambda)$.

THEOREM 4.6. *If the measure $\mu$ has a finite entropy, then a $\mu$-boundary is maximal if and only if the asymptotic entropy of almost all associated conditional random walks vanishes.*

This result generalizes the entropy criterion of the triviality of the Poisson boundary due to A. M. Vershik and the author [KV83] (announced in 1979) and to Y. Derriennic [De80]. In view of Theorem 4.6, in order to prove the maximality of a given $\mu$-boundary one has to show that, with probability bounded away from zero, the one-dimensional distributions of the conditional random walks are concentrated on subsets of subexponential growth. It leads to two simple purely geometric criteria of boundary maximality. Both require an approximation of the sample paths of the random walk in terms of their limit behaviour. For simplicity assume that $G$ is finitely generated, and denote by $d$ the left-invariant metric corresponding to a word length $\delta$ on $G$. Let $B$ be a $\mu$-boundary, and $\Pi$ be the corresponding projection from the path space onto $B$.

THEOREM 5.5 ("ray" or "unilateral" approximation). *If there is a family of measurable maps $\pi_n : B \to G$ such that $\mathbf{P}$-a.e. $d(x_n, \pi_n(\Pi \boldsymbol{x})) = o(n)$, then $B$ is maximal.*

The second criterion applies simultaneously to a $\mu$-boundary $B_+$ and to a $\check{\mu}$-boundary $B_-$ (where $\check{\mu}(g) = \mu(g)^{-1}$). Denote by $\mathcal{G}_n$ the balls of the word metric centered at the group identity.

662                           VADIM A. KAIMANOVICHTHEOREM 6.4 ("strip" or "bilateral" approximation). *If there exists a G-equivariant measurable map S assigning to pairs $(b_-, b_+) \in B_- \times B_+$ nonempty subsets ("strips") $S(b_-, b_+) \subset G$ such that for* a.e. $(b_-, b_+) \in B_- \times B_+$

$$(0.3) \qquad \frac{1}{n} \log \bigl| S(b_-, b_+) \cap \mathcal{G}_{\delta(x_n)} \bigr| \to 0$$

*in probability* **P** *with respect to* $\boldsymbol{x} = \{x_n\}$, *then both* $B_-$ *and* $B_+$ *are maximal.*

The "thinner" the strips $S(b_-, b_+)$, the larger the class of measures for which condition (0.3) is satisfied; i.e., sample paths $\{x_n\}$ may be allowed to go to infinity "faster". If the strips $S(\gamma_-, \gamma_+)$ grow *subexponentially* then condition (0.3) is satisfied for any probability measure $\mu$ with a *finite first moment*, and if the strips grow *polynomially* then (0.3) is satisfied for any measure $\mu$ with a finite *first logarithmic moment* $\sum \log \delta(g) \mu(g)$ (see Theorem 6.5).

As an application, we consider several classes of groups (loosely speaking, we call them "groups with hyperbolic properties"): Gromov word hyperbolic groups, groups with infinitely many ends, fundamental groups of rank 1 manifolds and discrete subgroups of semi-simple Lie groups. All these groups are endowed with natural boundaries (the hyperbolic boundary, the space of ends, the sphere at infinity, and the Furstenberg boundary of the ambient Lie group, respectively), and it is known that the sample paths of random walks on these groups converge to the natural boundaries. We show (Theorems 7.7, 8.4, 9.2, 10.7) that in fact the Poisson boundary for random walks on these groups can be identified with the natural boundaries under the condition that the measure $\mu$ has finite entropy and finite first logarithmic moment (in particular, if $\mu$ has a finite first moment). However, the problem of the identification of the Poisson boundary for an *arbitrary* measure on these groups still remains open.

The proofs are based on the fact that for all considered classes of groups there are natural "strips" of polynomial growth joining pairs of boundary points. These are geodesic pencils for hyperbolic groups, unions of cut sets separating two ends for groups with infinitely many ends, geodesics for rank 1 groups and geodesic flats for discrete subgroups of semi-simple Lie groups. Actually, the existence of strips can also be used for proving the boundary convergence (see Theorem 2.4). In combination with the strip criterion the latter gives general conditions (satisfied for word hyperbolic groups and for groups with infinitely many ends) which guarantee the "stochastic maximality" of a group compactification, i.e., that the sample paths of a random walk converge in this compactification and that the compactification boundary with the resulting hitting measure is isomorphic to the Poisson boundary of the random walk (see Theorem 6.6).

For the groups considered in the present paper existence of strips is almost evident, whereas checking the ray criterion (which, in a sense, amounts to proving an appropriate generalization of the Oseledec multiplicative ergodic theorem) may be rather elaborate if it does not fail altogether. For the sake of comparison we show how both criteria work for word hyperbolic groups (the ray approximation for Gromov hyperbolic spaces obtained in Theorem 7.2 is interesting on its own). Yet another application of the strip criterion is the identification of the Poisson boundary for the Teichmüller modular groups [KM96] (the natural boundary here is the Thurston boundary, and the strips in this case are the Teichmüller geodesics), where existence of a ray approximation remains an open question.

However, there are also situations when the ray criterion is more helpful than the strip criterion. One can see it already in the example of discrete subgroups of semi-simple Lie groups, where the strip criterion requires a non-degeneracy assumption on the measure $\mu$ which is not necessary at all for the ray approximation. A far reaching generalization of this example was recently obtained by Karlsson and Margulis [KM99] who proved the ray approximation for an arbitrary group of motions of a nonpositively curved space and used it for an identification of the Poisson boundary.

The paper has the following structure. Section 1 is devoted to background definitions and notations. In Section 2 the relationship between group compactifications and $\mu$-boundaries is discussed. In Section 3 we describe the conditional random walks with respect to a $\mu$-boundary, after which in Section 4 we obtain the entropy criterion of maximality of a $\mu$-boundary. The ray approximation and the strip approximation criteria are proved in Sections 5 and 6, respectively. In the remaining Sections 7–10 we consider the applications to concrete classes of groups.

A significant part of the paper was written during a stay at the University of Manchester. I would also like to thank the UNAM Institute of Mathematics at Cuernavaca, Mexico, where the paper was finished, for support and excellent working conditions.

## 1. Random walks on groups and the Poisson boundary

1.1. Let $G$ be a countable group, and $\mu$ be a probability measure on $G$. The (right) *random walk* on $G$ determined by the measure $\mu$ is the Markov chain on $G$ with the transition probabilities

$$(1.1) \qquad p(x,y) = \mu(x^{-1}y)$$

invariant with respect to the left action of the group $G$ on itself. Thus, the position $x_n$ of the random walk at time $n$ is obtained from its position $x_0$ at time 0 by multiplying by independent $\mu$-distributed right *increments* $h_i$:

$$(1.2) \qquad x_n = x_0 h_1 h_2 \cdots h_n .$$

The *Markov operator* $P = P_\mu$ of averaging with respect to the transition probabilities of the random walk $(G, \mu)$ is

$$Pf(x) = \sum_y p(x,y)f(y) = \sum_h \mu(h)f(xh) .$$

Its adjoint operator $\theta \mapsto \theta P$ acts on the space of measures $\theta$ on $G$ as the convolution with the measure $\mu$. If $\theta$ is the distribution of the position of the random walk at time $n$, then $\theta P = \theta \mu$ is the distribution of its position at time $n+1$.

*Here and below we use the notation $\alpha\beta$ to denote the convolution of a measure $\alpha$ on $G$ and a measure $\beta$ on a $G$-space $X$ (or, on the group $G$ itself).*

1.2. Denote by $G^{\mathbb{Z}_+}$ the space of *sample paths* $\boldsymbol{x} = \{x_n\}$, $n \geq 0$ endowed with the coordinate-wise action of $G$. *Cylinder* subsets of the path space are denoted

$$(1.3) \qquad C_{g_0,g_1,\ldots,g_n} = \{\boldsymbol{x} \in G^{\mathbb{Z}_+} : x_i = g_i,\ 0 \leq i \leq n\} = \bigcap_{i=0}^n C_{g_i}^i ,$$

where $C_g^i = \{\boldsymbol{x} \in G^{\mathbb{Z}_+} : x_i = g\}$ are the *one-dimensional cylinders*. An initial distribution $\theta$ on $G$ determines the *Markov measure* $\mathbf{P}_\theta$ in the path space which is the isomorphic image of the measure $\theta \otimes \bigotimes_{n=1}^\infty \mu$ under the map (1.2). The *one-dimensional distribution* of the measure $\mathbf{P}_\theta$ at time $n$ (i.e., its image under the projection $\boldsymbol{x} \mapsto x_n$) is $\theta P^n = \theta \mu_n$, where $\mu_n$ is the *n-fold convolution* of the measure $\mu$.

By $\mathbf{P}$ we denote the measure in the path space corresponding to the initial distribution concentrated at the group identity $e$ (this is the measure in the path space which is the most important for us). All measures $\mathbf{P}_\theta = \theta \mathbf{P}$ are dominated by the ($\sigma$-finite) measure $\mathbf{P}_m$, where $m$ is the counting measure on $G$. The space $(G^{\mathbb{Z}_+}, \mathbf{P}_m)$ is a *Lebesgue space*, which allows us to use in the sequel the standard ergodic theory technique of *measurable partitions* and *conditional measures* due to Rohlin (e.g., see [Ro67]).

1.3. Let $T : \{x_n\} \mapsto \{x_{n+1}\}$ be the *time shift* in the path space $G^{\mathbb{Z}_+}$. Then

$$(1.4) \qquad T\mathbf{P}_\theta = \mathbf{P}_{\theta P} = \mathbf{P}_{\theta \mu}$$

for any measure $\theta$ on $G$. In particular, the $\sigma$-finite measure $\mathbf{P}_m$ is $T$-invariant.

*Definition.* The space of ergodic components $\Gamma$ of the time shift $T$ in the path space $(G^{\mathbb{Z}_+}, \mathbf{P}_m)$ is called the *Poisson boundary* of the random walk $(G, \mu)$.

In a more detailed way, denote by $\mathcal{A}_T$ the $\sigma$-algebra of all measurable $T$-invariant sets (mod 0). Since $(G^{\mathbb{Z}_+}, \mathbf{P}_m)$ is a Lebesgue space, there is a (unique up to an isomorphism) measurable space $\Gamma$ (the *space of ergodic components*) and a map bnd : $G^{\mathbb{Z}_+} \to \Gamma$ such that the $\sigma$-algebra $\mathcal{A}_T$ coincides (mod 0) with the $\sigma$-algebra of bnd-preimages of measurable subsets of $\Gamma$. Denote by $\eta$ the corresponding measurable partition of the path space into the bnd-preimages of points from $\Gamma$. We shall call $\eta$ the *Poisson partition*. The coordinate-wise action of $G$ on the path space commutes with the shift $T$, hence it projects to a canonical $G$-action on $\Gamma$. The measure $\nu = \text{bnd}(\mathbf{P})$ on the Poisson boundary is called the *harmonic measure*.

*Below we shall always endow the Poisson boundary $\Gamma$ of the couple $(G, \mu)$ with the harmonic measure $\nu = \nu_e$ determined by the group identity $e$ as a starting point. Unless otherwise specified, no conditions are imposed either on the group $\text{gr}(\mu)$ or on the semigroup $\text{sgr}(\mu)$ generated by the support of the measure $\mu$.*

Since $\text{bnd}(\mathbf{P}) = \text{bnd}(T\mathbf{P})$ by the definition of $\Gamma$, (1.4) implies that

$$(1.5) \qquad \nu = \text{bnd}(T\mathbf{P}) = \text{bnd}(\mathbf{P}_\mu) = \mu\text{bnd}(\mathbf{P}) = \mu\nu = \sum_g \mu(g) g\nu \, ;$$

i.e., the harmonic measure $\nu$ is $\mu$-*stationary*. Therefore, the translation $g\nu$ is absolutely continuous with respect to $\nu$ for any $g \in \text{sgr}(\mu)$.

1.4. The Bernoulli shift in the space of increments of the random walk determines the measure preserving ergodic transformation

$$(1.6) \qquad (U\boldsymbol{x})_n = x_1^{-1} x_{n+1}$$

of the path space $(G^{\mathbb{Z}_+}, \mathbf{P})$. Since $T\boldsymbol{x} = x_1(U\boldsymbol{x})$, we have

LEMMA. *For $\mathbf{P}$-a.e. sample path $\boldsymbol{x} = \{x_n\} \in G^{\mathbb{Z}_+}$*

$$\text{bnd}\,\boldsymbol{x} = x_1 \text{bnd}\, U\boldsymbol{x} \, .$$

1.5. *Definition.* The quotient $(\Gamma_\xi, \nu_\xi)$ of the Poisson boundary $(\Gamma, \nu)$ with respect to a certain $G$-invariant measurable partition $\xi$ is called a $\mu$-*boundary*.

Another way of defining a $\mu$-boundary is to say that it is a $G$-space with a $\mu$-stationary measure $\lambda$ such that $x_n\lambda$ weakly converges to a $\delta$-measure for **P**-a.e. path $\{x_n\}$ of the random walk $(G,\mu)$ [Fu73]. We shall denote by $\text{bnd}_\xi$ the canonical projection

$$\text{bnd}_\xi : (G^{\mathbb{Z}_+}, \mathbf{P}) \to (\Gamma, \nu) \to (\Gamma_\xi, \nu_\xi) ,$$

and by $\eta_\xi$ the corresponding partition of the path space.

If $\Pi$ is a $T$-invariant equivariant measurable map from the path space $(G^{\mathbb{Z}_+}, \mathbf{P})$ to a $G$-space $B$, then $(B, \Pi\mathbf{P})$ is a $\mu$-boundary. For example, such a map arises in the situation when $G$ is embedded into a topological $G$-space $X$, and **P**-a.e. sample path $\boldsymbol{x} = \{x_n\}$ converges to a limit $\Pi\boldsymbol{x} \in X$.

1.6. *Definition.* A compactification of the group $G$ is called $\mu$-*maximal* if the sample paths of the random walk $(G,\mu)$ converge a.e. in this compactification, and the arising $\mu$-boundary is in fact isomorphic to the Poisson boundary of $(G,\mu)$.

This property means that the compactification is indeed *maximal in a measure theoretic sense*; i.e., there is no way (up to measure 0) of further splitting the boundary points of this compactification. Below we shall give general geometric criteria for maximality of $\mu$-boundaries and $\mu$-maximality of group compactifications using a quantitative approach based on the entropy theory of random walks.

## 2. Group compactifications and $\mu$-boundaries

2.1. Let $\overline{G} = G \cup \partial G$ be a compactification of a countable group $G$ which is compatible with the group structure on $G$ in the sense that the action of $G$ on itself by left translations extends to an action on $\overline{G}$ by homeomorphisms. We shall always assume that $\overline{G}$ is separable and introduce the following conditions on $\overline{G}$:

(CP) If a sequence $g_n \in G$ converges to a point from $\partial G$ in the compactification $\overline{G}$, then the sequence $g_n x$ converges to the same limit for any $x \in G$.

(CS) The boundary $\partial G$ consists of at least 3 points, and there is a $G$-equivariant Borel map $S$ assigning to pairs of distinct points $(b_1, b_2)$ from $\partial G$ nonempty subsets (*strips*) $S(b_1, b_2) \subset G$ such that for any 3 pairwise distinct points $\overline{b}_i \in \partial G$, $i = 0, 1, 2$ there exist neighbourhoods $\overline{b}_0 \in \mathcal{O}_0 \subset \overline{G}$ and $\overline{b}_i \in \mathcal{O}_i \subset \partial G$, $i = 1, 2$ with the property that

$$S(b_1, b_2) \cap \mathcal{O}_0 = \varnothing \qquad \text{for all } b_i \in \mathcal{O}_i, \ i = 1, 2 .$$

Condition (CP) is called *projectivity* in [Wo93], whereas condition (CS) means that points from $\partial G$ are *separated by the strips* $S(b_1, b_2)$. As we shall see below (see Theorem 6.6), it is often convenient to take for $S(b_1, b_2)$ the union of all bi-infinite geodesics in $G$ (provided with a Cayley graph structure) which have $b_1, b_2$ as their endpoints.

2.2. LEMMA. *Let $\overline{G} = G \cap \partial G$ be a compactification satisfying conditions* (CP), (CS), *and $(g_n) \subset G$ be a sequence such that $g_n \to \overline{b} \in \partial G$. Then for any nonatomic probability measure $\lambda$ on $\partial G$ the translations $g_n \lambda$ converge to the point measure $\delta_{\overline{b}}$ in the weak* topology.*

*Proof.* If $g_n b \to \overline{b}$ for all $b \in \partial G$, there is nothing to prove. Otherwise, passing to a subsequence we may assume that there exists $b_1 \in \partial G$ such that $g_n b_1 \to \overline{b}_1 \neq \overline{b}$. We claim that then $g_n b \to \overline{b}$ for all $b \neq b_1$. Indeed, if not, then passing again to a subsequence we may assume that there is $b_2 \neq b_1$ such that $g_n b_2 \to \overline{b}_2 \neq \overline{b}$. Take a point $x \in S(b_1, b_2)$, then by condition (CS) the only possible limit points of the sequence $g_n x$ are $\overline{b}_1$ or $\overline{b}_2$, which contradicts condition (CP). Since the measure $\lambda$ is nonatomic, the claim implies that $g_n \lambda \to \delta_{\overline{b}}$. Thus, any sequence $(g_n)$ with $g_n \to \overline{b}$ has a subsequence $(g_{n_k})$ with $g_{n_k} \lambda \to \delta_{\overline{b}}$, so that $g_n \lambda \to \delta_{\overline{b}}$. □

2.3. *Definition.* A subgroup $G' \subset G$ is called *elementary* with respect to a compactification $\overline{G} = G \cap \partial G$ if $G'$ fixes a finite subset of $\partial G$.

2.4. THEOREM. *Let $\overline{G} = G \cap \partial G$ be a separable compactification of a countable group $G$ satisfying conditions* (CP), (CS), *and $\mu$ be a probability measure on $G$ such that the subgroup $\mathrm{gr}\,(\mu)$ generated by its support is nonelementary with respect to this compactification. Then $\mathbf{P}$-a.e. sample path $\boldsymbol{x} = \{x_n\}$ converges to a limit $\Pi \boldsymbol{x} \in \partial G$. The limit measure $\lambda = \Pi \mathbf{P}$ is purely nonatomic, the measure space $(\partial G, \lambda)$ is a $\mu$-boundary, and $\lambda$ is the unique $\mu$-stationary probability measure on $\partial G$.*

*Proof.* By compactness of $\partial G$ there exists a $\mu$-stationary probability measure $\lambda$ on $\partial G$. The measure $\lambda$ is purely nonatomic. Indeed, let $m$ be the maximal weight of its atoms, and $A_m \subset \partial G$ be the finite set of atoms of weight $m$. Since $\lambda$ is $\mu$-stationary, $\lambda(b) = \sum_g \mu(g) \lambda(g^{-1} b)$ for any $b \in A_m$, whence $A_m$ is $\mathrm{sgr}\,(\mu)^{-1}$-invariant, which by finiteness of $A_m$ implies that $A_m$ is also $\mathrm{gr}\,(\mu)$-invariant, the latter being impossible because the group $\mathrm{gr}\,(\mu)$ is nonelementary.

Since the measure $\lambda$ is $\mu$-stationary, $\mathbf{P}$-a.e. sequence of measures $x_n \lambda$ converges weakly to a probability measure $\lambda(\boldsymbol{x})$ (see [Fu73], [Ma91, Chap. 6]). As $\mathrm{gr}\,(\mu)$ is nonelementary, a.e. sample path $\boldsymbol{x} = \{x_n\}$ is unbounded as a subset of $G$. Then by Lemma 2.2 $\Pi \boldsymbol{x} = \lim x_n$ a.e. exists, and $\lambda(\boldsymbol{x}) = \delta_{\Pi \boldsymbol{x}}$.

Put $\nu = \Pi \mathbf{P}$, so that $(\partial G, \nu)$ is a $\mu$-boundary (see §1.5). By $\mu$-stationarity of $\lambda$

$$\lambda = \mu_n \lambda = \sum_g \mu_n(g)\, g\lambda = \int x_n \lambda\, d\mathbf{P}(\boldsymbol{x}) \qquad \text{for all } n \geq 0\,.$$

Since $x_n \lambda \to \delta_{\Pi \boldsymbol{x}}$, by passing to the limit on $n$ we obtain that $\lambda = \nu$. □

COROLLARY. *Under conditions of Theorem 2.4 the action of $G$ on $\partial G$ is mean proximal.*

*Remark.* Our approach is different from the usual one which consists of deducing mean proximality from some contraction properties of the action (see [Fu73], [Ma91, Chap. 6], [Wo93] for the definition and examples).

## 3. Conditional measures and the Doob transform

3.1. Denote by $H_1^+(G,\mu)$ the convex set of all nonnegative harmonic functions $f$ on $\mathrm{sgr}\,(\mu)$ (i.e., such that $Pf = f$) normalized by the condition $f(e) = 1$. Any function $f \in H_1^+(G,\mu)$ determines a new Markov chain (the *Doob transform*) on $\mathrm{sgr}\,(\mu)$ whose transition probabilities

$$p^f(x,y) = \mu(x^{-1}y)\frac{f(y)}{f(x)}$$

are "cohomologous" to the transition probabilities (1.1) of the original random walk (e.g., see [Dy69] where the term "$h$-process" is used). Denote by $\mathbf{P}^f$ the associated Markov measure on $G^{\mathbb{Z}_+}$ (with the initial distribution $\delta_e$). For any cylinder subset (see equation (1.3))

$$\mathbf{P}^f(C_{e,g_1,\ldots,g_n}) = \mathbf{P}(C_{e,g_1,\ldots,g_n})f(g_n)\,,$$

so that the map $f \mapsto \mathbf{P}^f$ is affine.

3.2. If $A$ is a measurable subset of the Poisson boundary with $\nu(A) > 0$, then by the Markov property for any cylinder set $C_{e,g_1,\ldots,g_n}$

$$\mathbf{P}(C_{e,g_1,\ldots,g_n} \cap \mathrm{bnd}^{-1}A) = \mathbf{P}(C_{e,g_1,\ldots,g_n})\mathbf{P}_{g_n}(\mathrm{bnd}^{-1}A) = \mathbf{P}(C_{e,g_1,\ldots,g_n})g_n\nu(A)\,,$$

whence

$$\mathbf{P}(C_{e,g_1,\ldots,g_n}|\mathrm{bnd}^{-1}A) = \frac{\mathbf{P}(C_{e,g_1,\ldots,g_n})g_n\nu(A)}{\mathbf{P}(\mathrm{bnd}^{-1}A)} = \mathbf{P}(C_{e,g_1,\ldots,g_n})\frac{g_n\nu(A)}{\nu(A)}\,;$$

i.e., the conditional measure $\mathbf{P}^A(\cdot) = \mathbf{P}(\cdot|\mathrm{bnd}^{-1}A)$ is the Doob transform of the measure $\mathbf{P}$ determined by the normalized harmonic function $\varphi_A(x) = x\nu(A)/\nu(A)$. Now,

$$\varphi_A = \frac{1}{\nu(A)}\int_A \varphi_\gamma\, d\nu(\gamma)\,, \qquad \varphi_\gamma(x) = \frac{dx\nu}{d\nu}(\gamma)\,,$$

whence (by the fact that the Doob transform is affine)

$$\mathbf{P}^A = \frac{1}{\nu(A)} \int_A \mathbf{P}^\gamma \, d\nu(\gamma) \,,$$

where $\mathbf{P}^\gamma$ are the Doob transforms determined by the functions $\varphi_\gamma$ (their harmonicity follows from $\mu$-stationarity (see equation (1.5)) of the measure $\nu$), which by the definition of systems of conditional measures in Lebesgue spaces (e.g., see [Ro67]) means

THEOREM. *The measures*

$$\mathbf{P}^\gamma(C_{e,g_1,\ldots,g_n}) = \mathbf{P}(C_{e,g_1,\ldots,g_n}|\gamma) = \mathbf{P}(C_{e,g_1,\ldots,g_n}) \frac{dg_n\nu}{d\nu}(\gamma)$$

*corresponding to the Markov operators $P^\gamma$ on* sgr$(\mu)$ *with transition probabilities*

$$p^\gamma(x,y) = \mu(x^{-1}y) \frac{dy\nu}{dx\nu}(\gamma)$$

*are the canonical system of conditional measures of the measure $\mathbf{P}$ with respect to the Poisson boundary.*

3.3. Let now $(\Gamma_\xi, \nu_\xi)$ be a $\mu$-boundary. Then

$$(3.1) \quad \frac{dg\nu_\xi}{d\nu_\xi}(\gamma_\xi) = \int \frac{dg\nu}{d\nu}(\gamma) \, d\nu(\gamma|\gamma_\xi) \qquad \text{for all } g \in \text{sgr}(\mu), \ \nu_\xi\text{-a.e. } \gamma_\xi \in \Gamma_\xi \,,$$

where $\nu(\cdot|\gamma_\xi)$ are the conditional measures of the measure $\nu$ on the fibers of the projection $\Gamma \to \Gamma_\xi$, $\gamma \mapsto \gamma_\xi$. Denote by $\mathbf{P}^{\gamma_\xi} = \mathbf{P}^{\varphi_{\gamma_\xi}}$ the Doob transforms determined by the Radon-Nikodym derivatives $\varphi_{\gamma_\xi}(x) = dx\nu_\xi/d\nu_\xi(\gamma_\xi)$. Then (3.1) in combination with the fact that the Doob transform is affine implies that for $\nu_\xi$-a.e. $\gamma_\xi \in \Gamma_\xi$

$$\mathbf{P}^{\gamma_\xi} = \int \mathbf{P}^\gamma \, d\nu(\gamma|\gamma_\xi) \,,$$

whence by the transitivity of systems of conditional measures (e.g., see [Ro67]) we have

THEOREM. *The family of measures $\mathbf{P}^{\gamma_\xi}$, $\gamma_\xi \in \Gamma_\xi$ is the family of conditional measures of the measure $\mathbf{P}$ with respect to the $\mu$-boundary $(\Gamma_\xi, \nu_\xi)$.*

## 4. Entropy of conditional walks and maximality of $\mu$-boundaries

4.1. *From now on we shall assume that the measure $\mu$ has finite entropy $H(\mu) = -\sum_g \mu(g) \log \mu(g)$.* Then there exists the limit $h(G, \mu) = \lim_n H(\mu_n)/n$ which is called the *entropy* of the random walk $(G, \mu)$ (see [De80], [KV83]).

Denote by $\alpha_1^k$ the partition of the path space $(G^{\mathbb{Z}_+}, \mathbf{P})$ determined by the positions of the random walk at times $1, 2, \ldots, k$ (i.e., two sample paths $\boldsymbol{x}, \boldsymbol{x}'$ belong to the same class of $\alpha_1^k$ if and only if $x_i = x_i'$ for all $i = 1, 2, \ldots, k$), and put $\alpha = \alpha_1^1$. Given two measurable partitions $\xi, \zeta$ of the space $(G^{\mathbb{Z}_+}, \mathbf{P})$ we shall denote by $H_{\mathbf{P}}(\xi)$ (resp., $H_{\mathbf{P}}(\xi|\zeta)$) the *entropy* of $\xi$ (resp., the *conditional entropy* of $\xi$ with respect to $\zeta$); see [Ro67] for the definitions. Since the increments of the random walk are independent and $\mu$-distributed, $H_{\mathbf{P}}(\alpha_1^k) = kH(\mu) = kH_{\mathbf{P}}(\alpha)$.

4.2. Let $\xi$ be a $G$-invariant partition of the Poisson boundary, and $(\Gamma_\xi, \nu_\xi)$ be the corresponding $\mu$-boundary.

LEMMA. *For any $k \geq 1$*

$$H_{\mathbf{P}}(\alpha_1^k | \eta_\xi) = kH_{\mathbf{P}}(\alpha | \eta_\xi) = k\left[H(\mu) - \int \log \frac{dx_1 \nu_\xi}{d\nu_\xi}(\text{bnd}_\xi \boldsymbol{x}) \, d\mathbf{P}(\boldsymbol{x})\right].$$

*Proof.* Given a path $\boldsymbol{x} = \{x_n\} \in G^{\mathbb{Z}_+}$, the element of the partition $\alpha_1^k$ containing $\boldsymbol{x}$ is the cylinder $C_{e,x_1,\ldots,x_k}$, and the image of $\boldsymbol{x}$ in $\Gamma_\xi$ is $\text{bnd}_\xi \boldsymbol{x}$, whence by Theorem 3.3 the corresponding conditional probability is

$$\mathbf{P}(C_{e,x_1,\ldots,x_k} | \text{bnd}_\xi \boldsymbol{x}) = \mathbf{P}(C_{e,x_1,\ldots,x_k}) \frac{dx_k \nu_\xi}{d\nu_\xi}(\text{bnd}_\xi \boldsymbol{x}),$$

and

$$H_{\mathbf{P}}(\alpha_1^k | \eta_\xi) = -\int \log \mathbf{P}(C_{e,x_1,\ldots,x_k} | \text{bnd}_\xi \boldsymbol{x}) \, d\mathbf{P}(\boldsymbol{x})$$

$$= kH(\mu) - \int \log \frac{dx_k \nu_\xi}{d\nu_\xi}(\text{bnd}_\xi \boldsymbol{x}) \, d\mathbf{P}(\boldsymbol{x}).$$

Now, passing to the increments $h_n$ by (1.2), telescoping and using Lemma 1.4 we get

(4.1)
$$\frac{dx_k \nu_\xi}{d\nu_\xi}(\text{bnd}_\xi \boldsymbol{x}) = \frac{dh_1 \ldots h_k \nu_\xi}{d\nu_\xi}(\text{bnd}_\xi \boldsymbol{x})$$

$$= \prod_{i=1}^{k} \frac{dh_i \nu_\xi}{d\nu_\xi}(x_{i-1}^{-1} \text{bnd}_\xi \boldsymbol{x}) = \prod_{i=1}^{k} \frac{d(U^{i-1}\boldsymbol{x})_1 \nu_\xi}{d\nu_\xi}(\text{bnd}_\xi U^{i-1}\boldsymbol{x}).$$

Since the measure $\mathbf{P}$ is $U$-invariant, the claim follows. □

4.3. THEOREM. *Let $\xi \preccurlyeq \xi'$ be two $G$-invariant measurable partitions of the Poisson boundary $(\Gamma, \nu)$. Then $H_{\mathbf{P}}(\alpha | \eta_\xi) \geq H_{\mathbf{P}}(\alpha | \eta_{\xi'})$, and the equality holds if and only if $\xi = \xi'$.*

*Proof.* Obviously, if $\xi \preccurlyeq \xi'$, then $\eta_\xi \preccurlyeq \eta'_{\xi'}$, so that the inequality follows from [Ro67, Property 5.10]. If $H_{\mathbf{P}}(\alpha|\eta_\xi) = H_{\mathbf{P}}(\alpha|\eta_{\xi'})$, then by Lemma 4.2 $H_{\mathbf{P}}(\alpha_1^k|\eta_\xi) = H_{\mathbf{P}}(\alpha_1^k|\eta_{\xi'})$ for any $k \geq 1$. Therefore, again by [Ro67, Property 5.10], all finite dimensional distributions of the conditional measures $\mathbf{P}^{\gamma_\xi}$ and $\mathbf{P}^{\gamma_{\xi'}}$ coincide for $\nu$-a.e. point $\gamma \in \Gamma$; i.e., $\mathbf{P}^{\gamma_\xi} = \mathbf{P}^{\gamma_{\xi'}}$, which means that $\xi = \xi'$. □

4.4. *Definition.* A probability measure $\Lambda$ on $G^{\mathbb{Z}_+}$ has *asymptotic entropy* $\mathbf{h}(\Lambda)$ if it has the following Shannon-Breiman-McMillan type *equidistribution property*:

$$-\frac{1}{n}\log \Lambda(C_{x_n}^n) \to \mathbf{h}(\Lambda)$$

for $\Lambda$-a.e. $\boldsymbol{x} = \{x_n\} \in G^{\mathbb{Z}_+}$ and in the space $L^1(\Lambda)$.

4.5. THEOREM. *Let $\xi$ be a measurable $G$-invariant partition of the Poisson boundary $(\Gamma, \nu)$. Then for $\nu_\xi$-a.e. point $\gamma_\xi \in \Gamma_\xi$*

$$\mathbf{h}(\mathbf{P}^{\gamma_\xi}) = H_{\mathbf{P}}(\alpha|\eta_\xi) - H_{\mathbf{P}}(\alpha|\eta) \ .$$

*Proof.* We have to check that for $\nu_\xi$-a.e. point $\gamma_\xi \in \Gamma_\xi$

$$-\frac{1}{n}\log \mathbf{P}^{\gamma_\xi}(C_{x_n}^n) \to H_{\mathbf{P}}(\alpha|\eta_\xi) - H_{\mathbf{P}}(\alpha|\eta)$$

for $\mathbf{P}^{\gamma_\xi}$-a.e. sample path $\boldsymbol{x} = \{x_n\}$ and in the space $L^1(\mathbf{P}^{\gamma_\xi})$. Since $\mathbf{P}^{\gamma_\xi}$ are the conditional measures of the measure $\mathbf{P}$ with respect to $\Gamma_\xi$, it amounts to proving that

$$-\frac{1}{n}\log \mathbf{P}^{\mathrm{bnd}_\xi \boldsymbol{x}}(C_{x_n}^n) \to H_{\mathbf{P}}(\alpha|\eta_\xi) - H_{\mathbf{P}}(\alpha|\eta)$$

$\mathbf{P}$-a.e. and in the space $L^1(\mathbf{P})$. By Theorem 3.3

$$\frac{1}{n}\log \mathbf{P}^{\mathrm{bnd}_\xi \boldsymbol{x}}(C_{x_n}^n) = \frac{1}{n}\log \mathbf{P}(C_{x_n}^n) + \frac{1}{n}\log \frac{dx_n\nu_\xi}{d\nu_\xi}(\mathrm{bnd}_\xi \boldsymbol{x}) \ .$$

Since $\mathbf{h}(\mathbf{P}) = h(G, \mu)$ (see [De80], [KV83]), the first term in the right-hand side converges to $-h(G, \mu)$. On the other hand, telescoping as in (4.1), applying the Birkhoff Ergodic Theorem to the transformation $U$, and using Lemma 4.2, we obtain that the second term in the right-hand side converges to $H(\mu) - H_{\mathbf{P}}(\alpha|\eta_\xi)$. It remains to use the fact that $H_{\mathbf{P}}(\alpha|\eta) = H(\mu) - h(G, \mu)$ (see [KV83]). □

4.6. It is proved in [De80], [KV83] that $\mathbf{h}(\mathbf{P}) = 0$ if and only if the Poisson boundary is trivial. Combining Theorem 4.3 (where the point partition of the Poisson boundary is taken for $\xi'$) with Theorem 4.5 we get the following generalization of that criterion:

THEOREM. *A $\mu$-boundary $(B, \lambda) \cong (\Gamma_\xi, \nu_\xi)$ is the Poisson boundary if and only if $\mathbf{h}(\mathbf{P}^{\gamma_\xi}) = 0$ for almost all conditional measures of the measure $\mathbf{P}$ with respect to $\Gamma_\xi$.*

COROLLARY. *A $\mu$-boundary $(B, \lambda) \cong (\Gamma_\xi, \nu_\xi)$ is the Poisson boundary if and only if for $\nu_\xi$-a.e. point $\gamma_\xi \in \Gamma_\xi$ there exist $\varepsilon > 0$ and a sequence of sets $A_n = A_n(\gamma_\xi) \subset G$ such that*

(i) $\log |A_n| = o(n)$,
(ii) *$p_n^{\gamma_\xi}(A_n) > \varepsilon$ for all sufficiently large $n$, where $p_n^{\gamma_\xi}(g) = \mathbf{P}^{\gamma_\xi}(C_g^n)$ are the one-dimensional distributions of the measures $\mathbf{P}^{\gamma_\xi}$.*

## 5. Ray approximation

5.1. *Definition.* An increasing sequence $\mathcal{G} = (\mathcal{G}_k)_{k \geq 1}$ of sets exhausting a countable group $G$ is called a *gauge* on $G$. By

$$|g| = |g|_\mathcal{G} = \min\{k : g \in \mathcal{G}_k\}$$

we denote the corresponding *gauge function*. We shall say that a gauge $\mathcal{G}$ is

- *symmetric* if all gauge sets $\mathcal{G}_k$ are symmetric, i.e., $|g| = |g^{-1}|$ for all $g \in G$,
- *subadditive* if $|g_1 g_2| \leq |g_1| + |g_2|$ for all $g_1, g_2 \in G$,
- *finite* if all gauge sets are finite,
- *temperate* if it is finite and the gauge sets grow at most exponentially: $\sup_k \frac{1}{k} \log \operatorname{card} \mathcal{G}_k < \infty$.

A family of gauges $\mathcal{G}^\alpha$ is *uniformly temperate* if $\sup_{\alpha, k} \frac{1}{k} \log \operatorname{card} \mathcal{G}_k^\alpha < \infty$. Clearly, the family of *translations* $g\mathcal{G} = (g\mathcal{G}_k)$, $g \in G$ of any temperate gauge is uniformly temperate.

*The gauges considered below are not assumed to be finite or subadditive unless otherwise specified.*

An important class of gauges consists of *word gauges* (see [Gu80]), i.e., gauges $(\mathcal{G}_k)$ such that $\mathcal{G}_1$ is a set generating $G$ as a semigroup, and $\mathcal{G}_k = (\mathcal{G}_1)^k$ is the set of words of *length* $\leq k$ in the alphabet $\mathcal{G}_1$. Any word gauge is subadditive. It is symmetric if and only if the set $\mathcal{G}_1$ is symmetric, and finite if and only if $\mathcal{G}_1$ is finite. In the latter case the gauge is temperate. Any two finite word gauges $\mathcal{G}, \mathcal{G}'$ on a finitely generated group $G$ are *equivalent* (quasi-isometric) in the sense that there is a constant $C > 0$ such that

$$\frac{1}{C}|g|_{\mathcal{G}'} \leq |g|_\mathcal{G} \leq C|g|_{\mathcal{G}'} \qquad \text{for all } g \in G.$$

Thus, for a probability measure $\mu$ on a finitely generated group $G$ finiteness of its *first moment* $\sum_g |g|\mu(g)$, or of its *first logarithmic moment* $\sum_g \log|g|\mu(g)$ are independent of the choice of a finite word gauge $|\cdot|$ on $G$.

5.2. LEMMA (cf. [De86]) . *If $\mathcal{G}$ is a temperate gauge, and $|\mu|_\mathcal{G} < \infty$, then $H(\mu) < \infty$.*

5.3. For subadditive gauges the Kingman Subadditive Ergodic Theorem immediately implies (cf. [Gu80], [De80]):

LEMMA. *If $\mathcal{G}$ is a subadditive gauge, and $|\mu|_\mathcal{G} < \infty$, then the limit (rate of escape)*

$$\ell(G, \mu, \mathcal{G}) = \lim_{n \to \infty} \frac{|x_n|_\mathcal{G}}{n}$$

*exists for $\mathbf{P}$-a.e. sample path $\{x_n\}$ and in the space $L^1(\mathbf{P})$.*

5.4. THEOREM. *Let $\mu$ be a probability measure with finite entropy $H(\mu)$ on a countable group $G$, and $(B, \lambda) \cong (\Gamma_\xi, \nu_\xi)$ be a $\mu$-boundary. Denote by $\Pi = \mathrm{bnd}_\xi$ the projection from the path space $(G^{\mathbb{Z}_+}, \mathbf{P})$ to $(B, \lambda)$. If for $\lambda$-a.e. point $b \in B$ there exists a sequence of uniformly temperate gauges $\mathcal{G}^n = \mathcal{G}^n(b)$ such that*

(5.1) $$\frac{1}{n}|x_n|_{\mathcal{G}^n(\Pi \boldsymbol{x})} \to 0$$

*for $\mathbf{P}$-a.e. sample path $\boldsymbol{x} = \{x_n\}$, then $(B, \lambda)$ is the Poisson boundary of the pair $(G, \mu)$.*

*Proof.* Condition (5.1) is equivalent to saying that $|x_n|_{\mathcal{G}^n(b)}/n \to 0$ for $\lambda$-a.e. $b \in B$ and $\mathbf{P}^b$-a.e. sample path of the random walk conditioned by $b$ (see Theorem 3.3). Thus, $(B, \lambda)$ is the Poisson boundary by Theorem 4.6. $\square$

5.5. Now let $\pi_n : B \to G$ be a sequence of measurable maps from a $\mu$-boundary $B$ to the group $G$. Geometrically, one can think about the sequences $\pi_n(b)$, $b \in B$ as *"rays"* in $G$ corresponding to points from $B$. Taking in Theorem 5.4 $\mathcal{G}^n(b) = \pi_n(b)\mathcal{G}$, where $\mathcal{G}$ is a fixed temperate gauge on $G$, we obtain

THEOREM. *Let $\mu$ be a probability measure with finite entropy $H(\mu)$ on a countable group $G$, and $(B, \lambda) = \Pi(G^{\mathbb{Z}_+}, \mathbf{P})$ be a $\mu$-boundary. If there exist a temperate gauge $\mathcal{G}$ and a sequence of measurable maps $\pi_n : B \to G$ such that*

$$\frac{1}{n}\left|\left(\pi_n(\Pi \boldsymbol{x})\right)^{-1} x_n\right|_\mathcal{G} \to 0$$

*for $\mathbf{P}$-a.e. sample path $\boldsymbol{x} = \{x_n\}$, then $(B, \lambda)$ is the Poisson boundary of the pair $(G, \mu)$.*

## 6. Strip approximation

6.1. We have defined the path space $(G^{\mathbb{Z}_+}, \mathbf{P})$ (see §1.2) as the image of the space of independent $\mu$-distributed increments $\{h_n\}$, $n \geq 1$ under the map

$$\text{(6.1)} \qquad x_n = \begin{cases} e, & n = 0 \\ x_{n-1} h_n, & n \geq 1 . \end{cases}$$

Extending the relation $x_n = x_{n-1} h_n$ to all indices $n \in \mathbb{Z}$ (and always assuming that $x_0 = e$) we obtain the measure space $(G^{\mathbb{Z}}, \overline{\mathbf{P}})$ of *bilateral paths* $\overline{x} = \{x_n, n \in \mathbb{Z}\}$ corresponding to bilateral sequences of independent $\mu$-distributed increments $\{h_n\}$, $n \in \mathbb{Z}$. For negative indices $n$ formula (6.1) can be rewritten as

$$x_{-n} = x_{-n+1} h_{-n+1}^{-1}, \qquad n \geq 0 ,$$

so that

$$\check{x}_n = x_{-n} = h_0^{-1} h_{-1}^{-1} \cdots h_{-n+1}^{-1}, \qquad n \geq 0$$

is a sample path of the random walk on $G$ governed by the *reflected measure* $\check{\mu}(g) = \mu(g^{-1})$. The unilateral paths $x = \{x_n\}$, $n \geq 0$ and $\check{x} = \{\check{x}_n\} = \{x_{-n}\}$, $n \geq 0$ are independent; i.e., the map $\overline{x} \mapsto (x, \check{x})$ is an isomorphism of the measure spaces $(G^{\mathbb{Z}}, \overline{\mathbf{P}})$ and $(G^{\mathbb{Z}_+}, \mathbf{P}) \times (G^{\mathbb{Z}_+}, \check{\mathbf{P}})$, where $\check{\mathbf{P}}$ is the measure in the space of unilateral sample paths of the random walk $(G, \check{\mu})$.

6.2. Denote by $\overline{U}$ the measure preserving transformation of the space of bilateral paths $(G^{\mathbb{Z}}, \overline{\mathbf{P}})$ induced by the bilateral Bernoulli shift in the space of increments. It is the natural extension of the transformation $U$ of the unilateral path space $(G^{\mathbb{Z}_+}, \mathbf{P})$ defined in 1.4 and acts by the same formula (1.6) extended to all indices $n \in \mathbb{Z}$: for any $k \in \mathbb{Z}$

$$\text{(6.2)} \qquad (\overline{U}^k \overline{x})_n = x_k^{-1} x_{n+k} \qquad \text{for all } n \in \mathbb{Z} ;$$

i.e., the path $\overline{U}^k \overline{x}$ is obtained from the path $\overline{x}$ by translating it both in time (by $k$) and in space (by multiplying by $x_k^{-1}$ on the left in order to satisfy the condition $(\overline{U}^k \overline{x})_0 = e$). In terms of the unilateral paths $x$ and $\check{x}$ applying $\overline{U}^k$ consists (for $k > 0$) of canceling the first $k$ factors $x_k = h_1 h_2 \cdots h_k$ from the products $x_n = h_1 h_2 \cdots h_k \cdots h_n$, $n > 0$ (i.e., in applying to $x$ the transformation $U^k$) and adding on the left $k$ factors $x_k^{-1} = h_k^{-1} \cdots h_2^{-1} h_1^{-1}$ to the products $\check{x}_n = x_{-n} = h_0^{-1} h_{-1}^{-1} \cdots h_{-n+1}^{-1}$:

$$\overbrace{\cdots, h_{-1}, h_0, \underbrace{h_1, \cdots, h_{k-1}, h_k}}, \underbrace{h_{k+1}, \cdots} .$$

6.3. Denote by $\check{\Gamma}$ the Poisson boundary of the measure $\check{\mu}$, and by $\check{\nu}$ the corresponding harmonic measure.

THEOREM. *The action of the group $G$ on the product $(\check\Gamma \times \Gamma, \check\nu \otimes \nu)$ is ergodic.*

*Proof.* Denote by $\pi$ the measure preserving projection $\overline{\boldsymbol{x}} \mapsto (\check{\boldsymbol{x}}, \boldsymbol{x}) \mapsto (\operatorname{bnd}\check{\boldsymbol{x}}, \operatorname{bnd}\boldsymbol{x})$ from the bilateral path space $(G^{\mathbb{Z}}, \overline{\mathbf{P}})$ to the product space $(\check\Gamma \times \Gamma, \check\nu \otimes \nu)$. Then as it follows from formula (6.2), for any $k \in \mathbb{Z}$

$$(6.3) \qquad \pi(\overline{U}^k \overline{\boldsymbol{x}}) = x_k^{-1} \pi(\overline{\boldsymbol{x}})$$

(cf. Lemma 1.4). Now, if $A \subset \check\Gamma \times \Gamma$ is a $G$-invariant subset of $\check\Gamma \times \Gamma$ with $0 < \check\nu \otimes \nu(A) < 1$, then by (6.3) the preimage $\pi^{-1}(A)$ is $\overline{U}$-invariant with $0 < \overline{\mathbf{P}}(\pi^{-1}A) = \check\nu \otimes \nu(A) < 1$, which is impossible by the ergodicity of the bilateral Bernoulli shift $\overline{U}$. □

6.4. THEOREM. *Let $\mu$ be a probability measure with finite entropy $H(\mu)$ on a countable group $G$, and let $(B_-, \lambda_-)$ and $(B_+, \lambda_+)$ be $\check\mu$- and $\mu$-boundaries, respectively. If there exist a gauge $\mathcal{G} = (\mathcal{G}_k)$ on the group $G$ with gauge function $|\cdot| = |\cdot|_{\mathcal{G}}$ and a measurable $G$-equivariant map $S$ assigning to pairs of points $(b_-, b_+) \in B_- \times B_+$ nonempty "strips" $S(b_-, b_+) \subset G$ such that for all $g \in G$ and $\lambda_- \otimes \lambda_+$–a.e. $(b_-, b_+) \in B_- \times B_+$*

$$(6.4) \qquad \frac{1}{n} \log \operatorname{card} \left[ S(b_-, b_+) g \cap \mathcal{G}_{|x_n|} \right] \xrightarrow[n \to \infty]{} 0$$

*in probability with respect to the measure $\mathbf{P}$ in the space of sample paths $\boldsymbol{x} = \{x_n\}_{n \geq 0}$, then the boundary $(B_+, \lambda_+)$ is maximal.*

*Proof.* Denote by $\Pi_- : \overline{\boldsymbol{x}} \mapsto \check{\boldsymbol{x}} \mapsto \operatorname{bnd}_{\check\xi} \check{\boldsymbol{x}}$ and $\Pi_+ : \overline{\boldsymbol{x}} \mapsto \boldsymbol{x} \mapsto \operatorname{bnd}_\xi \boldsymbol{x}$ the projections of the bilateral path space $(G^{\mathbb{Z}}, \overline{\mathbf{P}})$ onto the boundaries $(B_-, \lambda_-) \cong (\check\Gamma_{\check\xi}, \check\nu_{\check\xi})$ and $(B_+, \lambda_+) \cong (\Gamma_\xi, \nu_\xi)$, respectively (cf. the proof of Theorem 6.3). Replacing if necessary the map $S$ with an appropriate right translation $(b_-, b_+) \mapsto S(b_-, b_+)g$, we may assume without loss of generality that

$$\lambda_- \otimes \lambda_+ \{(b_-, b_+) : e \in S(b_-, b_+)\} = \overline{\mathbf{P}}\left[e \in S(\Pi_-\overline{\boldsymbol{x}}, \Pi_+\overline{\boldsymbol{x}})\right] = p > 0 .$$

Using formula (6.3) in combination with the fact that the measure $\overline{\mathbf{P}}$ is $\overline{U}$-invariant, we then have for any $n \in \mathbb{Z}$

$$(6.5) \qquad \begin{aligned} \overline{\mathbf{P}}\left[x_n \in S(\Pi_-\overline{\boldsymbol{x}}, \Pi_+\overline{\boldsymbol{x}})\right] &= \overline{\mathbf{P}}\left[e \in x_n^{-1} S(\Pi_-\overline{\boldsymbol{x}}, \Pi_+\overline{\boldsymbol{x}})\right] \\ &= \overline{\mathbf{P}}\left[e \in S(x_n^{-1}\Pi_-\overline{\boldsymbol{x}}, x_n^{-1}\Pi_+\overline{\boldsymbol{x}})\right] \\ &= \overline{\mathbf{P}}\left[e \in S(\Pi_-\overline{U}^n\overline{\boldsymbol{x}}, \Pi_+\overline{U}^n\overline{\boldsymbol{x}})\right] \\ &= \overline{\mathbf{P}}\left[e \in S(\Pi_-\overline{\boldsymbol{x}}, \Pi_+\overline{\boldsymbol{x}})\right] = p . \end{aligned}$$

Since the image of the measure $\overline{\mathbf{P}}$ under the map $\overline{\boldsymbol{x}} \mapsto (\Pi_-\overline{\boldsymbol{x}}, \Pi_+\overline{\boldsymbol{x}})$ is $\lambda_- \otimes \lambda_+$, formula (6.5) can be rewritten as

$$(6.6) \qquad \int\int p_n^{b_+}\left[S(b_-, b_+)\right] d\lambda_-(b_-) d\lambda_+(b_+) = p ,$$

where $p_n^{b_+}$ are the one-dimensional distributions of the conditional measure $\mathbf{P}^{b_+}$.

Let
$$K_n = \min\{k \geq 1 : \mu_n(\mathcal{G}_k) \geq 1 - p/2\},$$

so that
$$\mathbf{P}\big[|x_n| \leq K_n\big] = \mu_n(\mathcal{G}_{K_n}) \geq 1 - p/2,$$

or, after conditioning by $\Pi_+ \overline{x}$,

(6.7) $$\int p_n^{b_+}\big[\mathcal{G}_{K_n}\big]\, d\lambda_+(b_+) \geq 1 - p/2.$$

Since for all $(b_-, b_+) \in B_- \times B_+$
$$p_n^{b_+}\big[S(b_-, b_+) \cap \mathcal{G}_{K_n}\big] \geq p_n^{b_+}\big[S(b_-, b_+)\big] + p_n^{b_+}\big[\mathcal{G}_{K_n}\big] - 1,$$

(6.6) and (6.7) imply
$$\int\int p_n^{b_+}\big[S(b_-, b_+) \cap \mathcal{G}_{K_n}\big]\, d\lambda_-(b_-) d\lambda_+(b_+) \geq p/2,$$

whence

(6.8) $$\lambda_- \otimes \lambda_+\Big\{(b_-, b_+) : p_n^{b_+}\big[S(b_-, b_+) \cap \mathcal{G}_{K_n}\big] \geq p/4\Big\} \geq p/4.$$

On the other hand, condition (6.4) implies that
$$\frac{1}{n} \log \operatorname{card}\big[S(b_-, b_+) \cap \mathcal{G}_{K_n}\big] \xrightarrow[n \to \infty]{} 0 \qquad \lambda_- \otimes \lambda_+ \text{-a.e.\ } (b_-, b_+) \in B_- \times B_+,$$

whence there exist a subset $Z \subset B_- \times B_+$ and a sequence $\varphi_n$ with $\log \varphi_n / n \to 0$ such that

(6.9) $$\lambda_- \otimes \lambda_+(Z) \geq 1 - p/8,$$

and

(6.10) $$\operatorname{card}\big[S(b_-, b_+) \cap \mathcal{G}_{K(n)}\big] \leq \varphi_n \qquad \text{for all } (b_-, b_+) \in Z.$$

Combining (6.8), (6.9) and (6.10) shows that there exists a sequence of sets $X_n \subset B_- \times B_+$ such that
$$\lambda_- \otimes \lambda_+(X_n) \geq p/8,$$
$$p_n^{b_+}\big[S(b_-, b_+) \cap \mathcal{G}_{K_n}\big] \geq p/4 \qquad \text{for all } (b_-, b_+) \in X_n,$$
$$\operatorname{card}\big[S(b_-, b_+) \cap \mathcal{G}_{K_n}\big] \leq \varphi_n \qquad \text{for all } (b_-, b_+) \in X_n.$$

Thus, taking $Y_n$ to be the projection of $X_n$ onto $B_+$, we have that $\lambda_+(Y_n) \geq p/8$, and for a.e. $b_+ \in Y_n$ there exists a set $A = A(b_+, n)$ with $p_n^{b_+}(A) \geq p/4$ and $\operatorname{card} A \leq \varphi_n$, so that the boundary $(B_+, \lambda_+)$ is maximal by Theorem 4.6. $\square$

6.5. Subexponentiality of the intersections $[S(b_-, b_+) \cap \mathcal{G}_{|x_n|}]$ is the key condition of Theorem 6.4. Thus, the "thinner" the strips $S(b_-, b_+)$, the larger the class of measures satisfying condition (6.4) of Theorem 6.4 (i.e., sample paths $\{x_n\}$ may be allowed to go to infinity "faster"). We shall illustrate this trade-off by giving two more operational corollaries to Theorem 6.4.

THEOREM. *Suppose that $\mathcal{G}$ is a subadditive temperate gauge on a countable group $G$ with gauge function $|\cdot| = |\cdot|_\mathcal{G}$ (particular case: $G$ is finitely generated, and $\mathcal{G}$ is a finite word gauge), and $\mu$ is a probability measure on $G$. Let $(B_-, \lambda_-)$ and $(B_+, \lambda_+)$ be $\check{\mu}$- and $\mu$-boundaries, respectively, and there exists a measurable $G$-equivariant map $B_- \times B_+ \ni (b_-, b_+) \mapsto S(b_-, b_+) \subset G$. If either*

(a) *the measure $\mu$ has a finite first moment $\sum |g| \mu(g)$, and for $\lambda_- \otimes \lambda_+$-a.e. $(b_-, b_+)$*

$$\frac{1}{k} \log \operatorname{card} \left[ S(b_-, b_+) \cap \mathcal{G}_k \right] \to 0$$

*(the strips $S(b_-, b_+)$ grow subexponentially with respect to $\mathcal{G}$), or*

(b) *the measure $\mu$ has a finite first logarithmic moment $\sum \log |g| \mu(g)$ and finite entropy $H(\mu)$, and for $\lambda_- \otimes \lambda_+$-a.e. $(b_-, b_+)$*

$$\sup_k \frac{1}{\log k} \log \operatorname{card} \left[ S(b_-, b_+) \cap \mathcal{G}_k \right] < \infty$$

*(the strips $S(b_-, b_+)$ grow polynomially),*

*then the boundaries $(B_-, \lambda_-)$ and $(B_+, \lambda_+)$ are maximal.*

*Proof.* (a) By Lemma 5.2, the measure $\mu$ has finite entropy, and by Lemma 5.3 there exists the rate of escape $\ell(G, \mu, \mathcal{G})$. Now, for any $g \in G$

$$\operatorname{card} \left[ S(b_-, b_+) g \cap \mathcal{G}_{|x_n|} \right] = \operatorname{card} \left[ S(b_-, b_+) \cap \mathcal{G}_{|x_n|} g^{-1} \right]$$
$$\leq \operatorname{card} \left[ S(b_-, b_+) \cap \mathcal{G}_{|x_n| + |g^{-1}|} \right],$$

whence condition (6.4) is satisfied.

(b) The proof is analogous to the proof of part (a), except that now we have to show that $\log |x_n|/n \to 0$. Indeed,

$$|x_n| = |h_1 h_2 \cdots h_n| \leq |h_1| + |h_2| + \cdots + |h_n|,$$

where $h_n$ are the independent $\mu$-distributed increments of the random walk. Since the measure $\mu$ has a finite first logarithmic moment, a.e. $\log |h_n|/n \to 0$, it follows that a.e. $\log |x_n|/n \to 0$. Now, for $\lambda_- \otimes \lambda_+$–a.e. $(b_-, b_+)$ and **P**-a.e. path $\{x_n\}$

$$\frac{1}{n} \log \operatorname{card} \left[ S(b_-, b_+) g \cap \mathcal{G}_{|x_n|} \right] \leq \frac{1}{n} \log \operatorname{card} \left[ S(b_-, b_+) \cap \mathcal{G}_{|x_n|+|g^{-1}|} \right]$$

$$= \frac{\log(|x_n| + |g^{-1}|)}{n} \cdot \frac{\log \operatorname{card} \left[ S(b_-, b_+) \cap \mathcal{G}_{|x_n|+|g^{-1}|} \right]}{\log(|x_n| + |g^{-1}|)} \to 0 \,. \quad \square$$

6.6. Let us introduce the following condition on a group compactification $\overline{G} = G \cup \partial G$.

(CG) There exists a left-invariant metric $d$ on $G$ such that the corresponding gauge $|\cdot|_d$ on $G$ is temperate and for any two distinct points $b_- \neq b_+ \in \partial G$

  (i) the pencil $P(b_-, b_+)$ of all $d$-geodesics $\alpha$ in $G$ such that $b_-$ (resp., $b_+$) is a limit point of the negative (resp., positive) ray of $\alpha$ is nonempty,

  (ii) there exists a finite set $A = A(b_-, b_+)$ such that any geodesic from the pencil $P(b_-, b_+)$ intersects $A(b_-, b_+)$.

Combining Theorems 2.4 and 6.5 then gives

THEOREM. *Let $\overline{G} = G \cap \partial G$ be a separable compactification of a countable group $G$ satisfying conditions* (CP), (CS), (CG), *and $\mu$ be a probability measure on $G$ such that*

 (i) *the subgroup $\operatorname{gr}(\mu)$ generated by its support is nonelementary with respect to this compactification,*

 (ii) *the measure $\mu$ has a finite entropy $H(\mu)$,*

(iii) *the measure $\mu$ has a finite first logarithmic moment with respect to the gauge determined by the metric $d$ from condition* (CG).

*Then the compactification $\overline{G}$ is $\mu$-maximal in the sense of Definition* 1.6.

*Proof.* Theorem 2.4 yields uniqueness of the measure $\lambda = \lambda_+$ and convergence, which implies that $(\partial G, \lambda_+)$ is a $\mu$-boundary. We shall deduce maximality of this boundary from Theorem 6.5. Indeed, the reflected measure $\check{\mu}$ satisfies the conditions of Theorem 2.4 simultaneously with the measure $\mu$. Let $\lambda_-$ be the unique $\check{\mu}$-stationary measure on $\partial G$. Since the measures $\lambda_-$ and $\lambda_+$ are purely nonatomic, the diagonal in $\partial G \times \partial G$ has zero measure $\lambda_- \otimes \lambda_+$, so that by condition (CG) for $\lambda_- \otimes \lambda_+$–a.e. $(b_-, b_+) \in \partial G \times \partial G$ there exists a minimal $M = M(b_-, b_+)$ such that all geodesics from the pencil $P(b_-, b_+)$ intersect a $M$-ball in $G$. Obviously, the map $(b_-, b_+) \mapsto M(b_-, b_+)$ is $G$-invariant, so that it must a.e. take a constant value $M_0$ by Theorem 6.3. Now define the strip $S(b_-, b_+) \subset G$ as the union of *all* balls $B$ of diameter $M_0$ such that any geodesic from the pencil $P(b_-, b_+)$ passes through $B$. This map is clearly $G$-equivariant, and for any geodesic $\alpha$ from $S(b_-, b_+)$ the strip $S(b_-, b_+)$ is contained in the $M_0$-neighbourhood of $\alpha$. Thus, the strips $S(b_-, b_+)$ have linear growth, so that the conditions of Theorem 6.5 are satisfied. □

## 7. Hyperbolic groups

7.1. In this section we identify the Poisson boundary for the Gromov hyperbolic groups (see [Gr87], [GH90] for the definitions). For the sake of comparison we shall use here both the ray and the strip approximations (Theorems 5.5 and 6.5, respectively).

*Definition.* A sequence of points $(x_n)$ in a Gromov hyperbolic space $X$ with metric $d$ is called *regular* if there exists a geodesic ray $\alpha$ and a number $l \geq 0$ (the *rate of escape*) such that $d(x_n, \alpha(nl)) = o(n)$, i.e., if the sequence $(x_n)$ asymptotically follows the ray $\alpha$. If $l > 0$, then we call $(x_n)$ a *nontrivial* regular sequence.

This notion is an analogue of the well-known notion of *Lyapunov regularity* (see [Ka89] and 10.2 below). The idea of the proof of the following result is due to T. Delzant. For Cartan-Hadamard manifolds with pinched sectional curvature another proof (using the Alexandrov Triangle Comparison Theorem) was given in [Ka85].

7.2. We shall fix a reference point $o \in X$, and put $|x| = d(o, x)$. Denote the Gromov product with respect to $o$ by $(x|y) = 1/2\big[|x| + |y| - d(x, y)\big]$.

THEOREM. *A sequence $(x_n)$ in a Gromov hyperbolic space $X$ is regular if and only if*

(i) $d(x_n, x_{n+1}) = o(n)$,
(ii) $|x_n|/n \to l \geq 0$.

*Proof.* Clearly, we just have to prove that (i) and (ii) imply regularity under the assumption that $l > 0$. Then $(x_{n-1}|x_n) = nl + o(n)$, and applying the quasi-metric $\rho(x, y) = \exp(-(x|y))$ (see [GH90]) yields convergence of $x_n$ to a point $x_\infty \in \partial X$ in the hyperbolic compactification of the space $X$. Now we fix geodesics $\alpha_n$ (resp., $\alpha_\infty$) joining the origin $o$ with the points $x_n$ (resp., $x_\infty$), and denote the points on these geodesics at distance $t$ from the origin by $[x_n]_t = \alpha_n(t)$ (resp., $[x_\infty]_t = \alpha_\infty(t)$).

Choose a positive number $\varepsilon < l/2$, and let

$$N = N(\varepsilon) = \min\{n > 0 : (x_{n-1}|x_n) \geq (l-\varepsilon)n\} \ .$$

In particular, $|x_n| \geq (l-\varepsilon)n$ for $n \geq N$, and the truncations $x_n^\varepsilon = [x_n]_{(l-\varepsilon)n}$ are well defined. The points $x_{n-1}^\varepsilon, x_n^\varepsilon$ belong to the sides of the geodesic triangle with vertices $o, x_{n-1}, x_n$, so that

$$d(x_{n-1}^\varepsilon, x_n^\varepsilon) \leq \big||x_{n-1}^\varepsilon| - |x_n^\varepsilon|\big| + 4\delta = l - \varepsilon + 4\delta \qquad \text{for all } n \geq N$$

because $|x_{n-1}^\varepsilon|, |x_n^\varepsilon| \leq (x_{n-1}|x_n)$, and geodesic triangles in $X$ are $4\delta$-thin [GH90, pp. 38, 41]. Therefore, for any two indices $n, m \geq N$

$$d(x_n^\varepsilon, x_m^\varepsilon) \leq |n - m|(l + 4\delta) \ ,$$
$$d(x_n^\varepsilon, x_m^\varepsilon) \geq \big||x_n^\varepsilon| - |x_m^\varepsilon|\big| = |n - m|(l - \varepsilon) \geq |n - m| l/2 \ .$$

It follows that the sequence $(x_n^\varepsilon)_{n \geq N}$ is a quasigeodesic, and by [GH90, p. 101] there exists a geodesic ray $\beta$ starting at the point $x_N^\varepsilon$ such that $d(x_n^\varepsilon, \beta) \leq H$ for any $n \geq N$ and a constant $H = H(\delta, l)$. Since $(x_n|x_n^\varepsilon) = n(l - \varepsilon) \to \infty$, the sequence $(x_n^\varepsilon)$ also converges to the point $x_\infty$, so that the geodesic rays $\beta$ and $\alpha_\infty$ are asymptotic. Thus, $d(x_n^\varepsilon, \alpha_\infty) \leq H + 8\delta$ and

$$d(x_n, \alpha_\infty) \leq H + 8\delta + d(x_n, x_n^\varepsilon) = H + 8\delta + \big(|x_n| - n(l - \varepsilon)\big)$$

for all sufficiently large $n$ (see [GH90, p. 117]). Since $\varepsilon$ can be made arbitrarily small, the claim is proven. □

7.3. A subgroup of a hyperbolic group $G$ is elementary with respect to the hyperbolic compactification in the sense of Definition 2.3 if and only if it is either finite or a finite extension of $\mathbb{Z}$; otherwise it is nonamenable (see [Gr87], [GH90]). Therefore, the Poisson boundary of a measure $\mu$ on $G$ is nontrivial if and only if gr$(\mu)$ is nonelementary, e.g. see [KV83]. Now assume that $\mu$ has a finite first moment and denote by $\ell$ the corresponding rate of escape (Lemma 5.3) with respect to a certain fixed word gauge $|\cdot|$ on $G$. If the Poisson boundary of $\mu$ is nontrivial, i.e., if gr$(\mu)$ is nonelementary, then by Theorem 5.5 $\ell > 0$. Furthermore, since the measure $\mu$ (i.e, the lengths of the increments $|h_n| = |x_{n-1}^{-1} x_n| = d(x_{n-1}, x_n)$) has finite first moment, $d(x_{n-1}, x_n) = o(n)$. Applying Theorem 7.2 we obtain

THEOREM. *Let $\mu$ be a probability measure with a finite first moment on a hyperbolic group $G$ such that the group* gr$(\mu)$ *is nonelementary. Then a.e. sample path of the random walk $(G, \mu)$ is a nontrivial regular sequence in $G$.*

7.4. For any $\xi \in \partial G$ choose a geodesic ray $\alpha_\xi$ from $e$ to $\xi$ in such a way that the map $\xi \mapsto \alpha_\xi$ is measurable (e.g., take for $\alpha_\xi$ the lexicographically minimal ray among all the rays joining $e$ and $\xi$), and let $\pi_n(\xi) = \alpha_\xi([n\ell])$, where $\ell$ is the rate of escape of the random walk $(G, \mu)$ and $[t]$ is the integer part of a number $t$. Then by Theorem 7.3 $d(x_n, \pi_n(x_\infty)) = o(n)$ for **P**-a.e. sample path $\{x_n\}$, so that by Theorem 5.5 we obtain

THEOREM. *Let $\mu$ be a probability measure with a finite first moment on a hyperbolic group $G$ such that the group $\mathrm{gr}\,(\mu)$ is nonelementary. Then a.e. sample path of the random walk $(G, \mu)$ converges in the hyperbolic compactification, and the hyperbolic boundary $\partial G$ with the resulting limit measure is isomorphic to the Poisson boundary of $(G, \mu)$.*

7.5. Using the strip approximation instead of the ray approximation allows us to obtain a stronger result in a simpler way.

PROPOSITION. *The hyperbolic compactification of a nonelementary hyperbolic group satisfies conditions* (CP), (CS), (CG) *from* 2.1 *and* 6.6.

*Proof.* Condition (CP) follows immediately from the definition of the hyperbolic compactification. For any two distinct points $\xi_- \neq \xi_+ \in \partial G$ let $S(\xi_-, \xi_+)$ be the union of points from all geodesics in $G$ joining $\xi_-$ and $\xi_+$. Then condition (CS) is implied by the quasi-convexity of geodesic hulls of subsets in the hyperbolic boundary [Gr87, 7.5.A], and condition (CG) follows from the fact that any two geodesics in a hyperbolic space with the same endpoints are within a uniformly bounded distance from each other. □

Applying Theorems 2.4 and 6.6 we then get

7.6. THEOREM. *Let $\mu$ be a probability measure on a hyperbolic group $G$ such that the subgroup $\mathrm{gr}\,(\mu)$ generated by its support is nonelementary. Then almost all sample paths $\{x_n\}$ converge to a (random) point $x_\infty \in \partial G$, so that $\partial G$ with the resulting limit measure $\lambda$ is a $\mu$-boundary. The measure $\lambda$ is the unique $\mu$-stationary probability measure on $\partial G$.*

7.7. THEOREM. *Under the conditions of Theorem 7.6, if the measure $\mu$ has finite entropy $H(\mu)$ and finite first logarithmic moment $\sum \mu(g) \log |g|$ (particular case: $\mu$ has a finite first moment), then $(\partial G, \lambda)$ is isomorphic to the Poisson boundary of $(G, \mu)$.*

*Remarks.* 1. Another proof of Theorem 7.6 in the case when $\mathrm{gr}\,(\mu) = G$ is given in [Wo93]. Unlike ours, it uses the contractivity of the $G$-action on $\partial G$.

2. If the measure $\mu$ is finitely supported and $\mathrm{sgr}\,(\mu) = G$, then the *Martin boundary* of the random walk coincides with the hyperbolic boundary [An90], so that in this case Theorems 7.6, 7.7 follow from the general Martin theory.

3. Theorems 7.6, 7.7 are easily seen to hold for any discrete discontinuous group of isometries of a Gromov hyperbolic space $X$ (in this case instead of a word gauge on $G$ one should take the gauge induced by the ambient metric on $X$).

## 8. Groups with infinitely many ends

8.1. The space of ends $\mathcal{E}(G)$ of a finitely generated group $G$ is defined as the space of ends of its Cayley graph with respect to a certain finite generating set $A$. Neither the space $\mathcal{E}(G)$ nor the end compactification $\overline{G} = G \cup \mathcal{E}(G)$ depend on the choice of $A$. For an end $\omega \in \mathcal{E}(G)$ and a finite set $K \subset G$ denote by $C(\omega, K)$ the connected component of $\overline{G} \setminus K$ containing $\omega$. The sets $C(\omega, K)$ form a basis of the end topology in $\overline{G}$ at the point $\omega$ (e.g., see [St71]). Clearly, any geodesic ray in $G$ converges to an end. Conversely, by standard compactness considerations for any two distinct ends from $\mathcal{E}(G)$ there exists a geodesic (not necessarily unique!) joining these ends.

The simplest example of a group with infinitely many ends is the free group $F_d$ of rank $d \geq 2$. This group is also hyperbolic. However, in general, a hyperbolic group may have a trivial space of ends (e.g., the fundamental group of a compact negatively curved manifold), and a group with infinitely many ends need not be hyperbolic (e.g., the free product of two copies of the group $\mathbb{Z}^2$). Nevertheless, groups with infinitely many ends still share with hyperbolic groups a number of geometric properties which are important for us.

8.2. LEMMA. *The end compactification of a finitely generated group with infinitely many ends satisfies conditions* (CP), (CS), (CG) *from* 2.1 *and* 6.6.

*Proof.* Condition (CP) is trivial. For verifying condition (CS) let $S(\omega_1, \omega_2)$ be the union of all geodesics in $G$ with endpoints $\omega_1 \neq \omega_2 \in \mathcal{E}(G)$. Take $\omega_0 \neq \omega_1, \omega_2 \in \mathcal{E}(G)$, then there is a finite set $K \subset G$ such that $C(\omega_0, K) \neq C(\omega_1, K), C(\omega_2, K)$, so that the intersection with $C(\omega_0, K)$ of any geodesic joining points in $C(\omega_1, K)$ and $C(\omega_2, K)$ must be contained in the (finite) union of all geodesic segments with endpoints from $K$. Finally, (CG) immediately follows from the definition of the space of ends. □

Now Theorems 2.4 and 6.6 imply

8.3. THEOREM. *Let $G$ be a finitely generated group with infinitely many ends, and $\mu$ be a probability measure such that the subgroup $\operatorname{gr}(\mu)$ generated by its support is nonelementary. Then almost all sample paths $\{x_n\}$ of the random walk $(G, \mu)$ converge to a (random) end $x_\infty \in \mathcal{E}(G)$, so that the space of ends $\mathcal{E}(G)$ with the resulting limit measure $\lambda$ is a $\mu$-boundary. The measure $\lambda$ is the unique $\mu$-stationary probability measure on $\mathcal{E}(G)$.*

8.4. THEOREM. *Under the conditions of Theorem 8.3, if the measure $\mu$ has finite entropy and finite first logarithmic moment (in particular, if $\mu$ has finite first moment), then the space $(\mathcal{E}(G), \lambda)$ is isomorphic to the Poisson boundary of the pair $(G, \mu)$.*

*Remark.* Our proof of Theorems 8.3, 8.4 is synthetic and does not evoke at all Stallings' structure theory of groups with infinitely many ends. If $\operatorname{gr}(\mu) = G$ Theorem 8.3 was proved by Woess [Wo93] using contractivity properties of the action of $G$ on $\mathcal{E}(G)$. A particular case of Theorem 8.4 when the measure $\mu$ is finitely supported and $\operatorname{sgr}(\mu) = G$ was proved by Woess [Wo89] by applying the Martin theory methods.

## 9. Fundamental groups of rank 1 manifolds

9.1. Let $M$ be a compact Riemannian manifold with nonpositive sectional curvature, and $\widetilde{M}$ be its universal covering space. Denote by $\partial \widetilde{M}$ the *sphere at infinity* of $\widetilde{M}$, which is the boundary of the *visibility compactification* of $\widetilde{M}$, e.g., see [Ba95]. The embedding $g \mapsto go \in \widetilde{M}$ (where $o \in \widetilde{M}$ is a fixed reference point) allows one to consider the visibility compactification of $\widetilde{M}$ as a compactification of the fundamental group $\pi_1(M)$. If $\widetilde{M}$ is *irreducible* (i.e., is not a product of two Cartan-Hadamard manifolds), then by the Rank Rigidity Theorem (e.g., see [Ba95]) $\widetilde{M}$ is either a symmetric space of noncompact type with rank at least 2, or $\widetilde{M}$ has a *regular* geodesic $\sigma$, i.e., such that there is no nontrivial parallel Jacobi field along $\sigma$ perpendicular to $\dot{\sigma}$. In the latter case $M$ is said to have *rank* 1. Note that the sectional curvature of a rank 1 manifold $M$ is not necessarily bounded away from 0, and its universal covering space $\widetilde{M}$ is not necessarily hyperbolic in the sense of Gromov.

9.2. THEOREM. *Let $\mu$ be a probability measure on the fundamental group $G = \pi_1(M)$ of a compact rank 1 Riemannian manifold $M$ such that $\operatorname{sgr}(\mu) = G$, and $\mu$ has a finite first logarithmic moment and finite entropy. Then a.e. sample path of the random walk $(G, \mu)$ converges in the visibility compactification, and the sphere $\partial \widetilde{M}$ with the resulting limit measure $\lambda$ is isomorphic to the Poisson boundary of the pair $(G, \mu)$.*

*Proof.* Convergence of sample paths was established by Ballmann for an arbitrary probability measure on $G$ with $\operatorname{sgr}(\mu) = G$ (see [Ba89, Theorem 2.2]). Moreover, $g_n \lambda \to \delta_\gamma$ weakly in $\widetilde{M} \cup \partial \widetilde{M}$ for any sequence $(g_n) \subset G$ such that $g_n \to \gamma$, i.e., the Dirichlet problem for $\mu$-harmonic functions with boundary data on $\partial \widetilde{M}$ is solvable (see [Ba89, Theorem 1.8]).

It remains to prove maximality of the $\mu$-boundary $(\partial \widetilde{M}, \lambda)$. Denote by $\lambda_+ = \lambda$ and $\lambda_-$ the hitting measures on $\partial \widetilde{M}$ determined by the measures $\mu$ and $\check{\mu}$, respectively. Since $M$ has rank one, the set $\mathcal{R} \subset \partial \widetilde{M} \times \partial \widetilde{M}$ of pairs of endpoints of regular bi-infinite geodesics in $\widetilde{M}$ is open nonempty, and for any pair of points $(\xi_-, \xi_+) \in \mathcal{R}$ there is a unique geodesic $\sigma(\xi_-, \xi_+)$ joining these points [Ba95]. Then solvability of the Dirichlet problem and quasi-invariance of the measures $\lambda_-, \lambda_+$ with respect to the action of $G$ (see 1.3) implies that $\lambda_- \times \lambda_+(\mathcal{R}) > 0$, whence $\lambda_- \times \lambda_+(\mathcal{R}) = 1$ by Theorem 6.3.

Since the quotient manifold $M$ is compact, there exists a number $d > 0$ such that for any point $x \in \widetilde{X}$ the $d$-ball centered at $x$ intersects the orbit $Go$. Then the strips in $G$ defined as

$$S(\xi_-, \xi_+) = \{g \in G : \operatorname{dist}(go, \sigma(\xi_-, \xi_+)) \leq d\}$$

are nonempty, and the map $(\xi_-, \xi_+) \mapsto S(\xi_-, \xi_+)$ is $G$-equivariant (here dist is the Riemannian metric on $\widetilde{M}$). The gauge $|g| = \operatorname{dist}(o, go)$ on $G$ is temperate and subadditive, and since $M$ is compact, it is equivalent to any finite word gauge on $G$, so that the measure $\mu$ has finite first logarithmic moment with respect to $|\cdot|$. Clearly, all strips $S(\xi_-, \xi_+)$ (being neighbourhoods of geodesics) have linear growth with respect to the gauge $|\cdot|$, and the conditions of Theorem 6.5 (b) are satisfied. □

*Remarks.* 1. For measures $\mu$ with a finite first moment Theorem 9.2 was first proved by Ballmann and Ledrappier [BL94]. Our "strip approximation" criterion (Theorems 6.4, 6.5) was inspired by the use of bilateral geodesics in [BL94], and the first part of our proof of Theorem 9.2 (existence of bilateral geodesics) is the same as in [BL94]. However, Theorem 6.5 allows us to obtain the result in greater generality and to avoid at the same time tedious dimension estimates (§3 in [BL94]).

2. In view of the rank rigidity theorem (e.g., see [Ba95]) Theorem 9.2 in combination with the results from Section 10 allows one to identify the Poisson boundary with (a subset of) the sphere at infinity for an arbitrary cocompact group of isometries of a Cartan-Hadamard manifold. Karlsson and Margulis [KM99] have recently obtained a generalization of the Oseledec multiplicative ergodic theorem which (with the help of the ray approximation criterion) leads to the same result for random walks with a finite first moment on an arbitrary discrete group of isometries of a Cartan-Hadamard manifold.

## 10. Discrete subgroups of semi-simple Lie groups

10.1. Let $\mathcal{G}$ be a connected semi-simple real Lie group with finite center, $\mathcal{K}$ be its maximal compact subgroup, and $S = \mathcal{G}/\mathcal{K}$ be the corresponding Riemannian symmetric space with the origin $o \cong \mathcal{K}$. Fix a dominant Weyl chamber $\mathfrak{A}^+$ in the Lie algebra $\mathfrak{A}$ of a Cartan subgroup $\mathcal{A}$, and denote by $\mathfrak{A}_1^+$ (resp., by $\overline{\mathfrak{A}}_1^+$) the intersection of $\mathfrak{A}^+$ (resp., of its closure $\overline{\mathfrak{A}}^+$) with the unit sphere of the Euclidean distance $\|\cdot\|$ determined by the Killing form $\langle \cdot, \cdot \rangle$. Any point $x \in S$ can be presented as $x = k(\exp a)o$, where $k \in \mathcal{K}$, and $a = r(x) \in \overline{\mathfrak{A}}^+$ is the uniquely determined *radial part* of $x$. Then the Riemannian distance $\operatorname{dist}(o, x)$ from $o$ to $x$ equals $\|r(x)\|$.

Denote by $\partial S$ the boundary (the *sphere at infinity*) of the *visibility compactification* of $S$ (cf. §9.1). We identify points from $\partial S$ with geodesic rays originating from $o$. The $\mathcal{G}$-orbits in $\partial S$ are parameterized by vectors $a \in \overline{\mathfrak{A}}_1^+$: the orbit $\partial S_a$ consists of the limits of all geodesic rays of the form $\xi(t) = g \exp(ta)o$. Stabilizers of points $\xi \in \partial S$ are parabolic subgroups of $\mathcal{G}$, which are minimal if and only if $\xi \in \partial S_a, a \in \mathfrak{A}_1^+$. Thus, the orbits $\partial S_a$ corresponding to nondegenerate vectors $a \in \mathfrak{A}_1^+$ are isomorphic to the *Furstenberg boundary* $\mathcal{B} = \mathcal{G}/\mathcal{P}$, where $\mathcal{P} = \mathcal{MAN}$ is the minimal parabolic subgroup determined by the Iwasawa decomposition $\mathcal{G} = \mathcal{KAN}$ (i.e., $\mathcal{M}$ is the centralizer of $\mathcal{A}$ in $\mathcal{K}$), and the orbits $\partial S_a$ corresponding to vectors $a$ from the walls of the Weyl chamber $\mathfrak{A}^+$ are isomorphic to quotients of the Furstenberg boundary (i.e., to quotients of $\mathcal{G}$ by nonminimal parabolic subgroups) (see [Ka89]). Moreover, there exists a canonical map $\overline{\mathfrak{A}}_1^+ \times \mathcal{B} \to \partial S$ such that $\{a\} \times \mathcal{B} \to \partial S_a$ is one-to-one for $a \in \mathfrak{A}_1^+$ (cf. below §10.4).

10.2. We call a sequence of points $x_n \in S$ *regular* if there exists a geodesic ray $\xi$ and a number $l \geq 0$ such that $\operatorname{dist}(x_n, \xi(nl)) = o(n)$ (cf. Definition 7.1). If $l > 0$, then $x_n$ converges in the visibility compactification to the same point as the ray $\xi$.

THEOREM ([Ka89]). *A sequence of points $x_n$ in a noncompact Riemannian symmetric space $S$ is regular if and only if $\operatorname{dist}(x_n, x_{n+1}) = o(n)$ and there exists a limit $a = \lim r(x_n)/n \in \overline{\mathfrak{A}}^+$.*

*Remark.* The definition of regular sequences is inspired by the notion of *Lyapunov regularity* (for the symmetric space $\operatorname{SL}(n, R)/\operatorname{SO}(n)$ these notions coincide, see §10.9), and Theorem 10.2 is a geometric counterpart of the *Oseledec multiplicative ergodic theorem* (see [Ka89]). Therefore we call the vector $a$ the *Lyapunov vector* of the sequence $x_n$.

10.3. THEOREM. *Let $\mu$ be a probability measure on a discrete subgroup $G \subset \mathcal{G}$ of a semi-simple Lie group $\mathcal{G}$ with a finite first moment*

$$\sum \text{dist}(o, go)\mu(g) < \infty.$$

*Then*

(i) **P**-*a.e. sample path $\{x_n\}$ of the random walk $(G, \mu)$ is regular, and the Lyapunov vector $a = a(\mu) = \lim r(x_n o)/n \in \overline{\mathfrak{A}}^+$ does not depend on $\{x_n\}$,*

(ii) *If $a \neq 0$, then for **P**-a.e. sample path $\{x_n\}$ the sequence $x_n o$ converges in the visibility compactification to a limit point from the orbit $\partial S_a$,*

(iii) *If $a = 0$, then the Poisson boundary of the pair $(G, \mu)$ is trivial, and if $a \neq 0$ it is isomorphic to $\partial S_a$ with the limit measure determined by* (ii).

*Proof.* Existence of the Lyapunov vector follows from Lemma 5.3 applied to matrix norms of finite dimensional representations of $\mathcal{G}$, see [Ka89]. Moreover, finiteness of the first moment of the measure $\mu$ implies that **P**-a.e. $\text{dist}(x_n o, x_{n+1} o) = o(n)$ (cf. the proof of Theorem 7.3), so that (i) and (ii) follow from Theorem 10.2.

Since the growth of $S$ is exponential, the gauge $g \mapsto \text{dist}(o, go)$ on $G$ induced by the Riemannian metric dist is temperate (see Definition 5.1), and combining Lemma 5.2 and Theorems 5.5, 10.2 we get (iii). □

*Remark.* If the group $\text{gr}(\mu)$ generated by the support of $\mu$ is nonamenable, then the Poisson boundary of $(G, \mu)$ is nontrivial (e.g., see [KV83]), and thereby $a \neq 0$.

10.4. If the measure $\mu$ does not have a finite first moment and the rank of $\mathcal{G}$ is greater than 1, convergence in the visibility compactification does not necessarily hold. However, in this situation one can use another compactification of the associated symmetric space by imposing some irreducibility conditions on the group $\text{gr}(\mu)$.

The map $go \mapsto gm$, where $m$ is the unique $\mathcal{K}$-invariant probability measure on $\mathcal{B}$, determines an embedding of $S$ into the space of Borel probability measures on $\mathcal{B}$, which gives rise to the *Satake-Furstenberg compactification* of $S$ obtained as the closure of the family of measures $\{gm\}$ in the weak topology. The boundary of this compactification consists (unless $S$ has rank 1) of several $G$-transitive components, one of which is $\mathcal{B}$ (corresponding to limit $\delta$-measures). If a sequence $x_n \in S$ converges in the visibility compactification to a point $b$ from a nondegenerate orbit $\partial S_a \cong \mathcal{B}$, $a \in \mathfrak{A}_1^+$, then $x_n$ also converges to $b \in \mathcal{B}$ in the Satake-Furstenberg compactification (see [Mo64]).

Another definition of the Furstenberg boundary $\mathcal{B}$ (analogous to that of the visibility boundary $\partial S$) can be given in terms of maximal totally geodesic flat subspaces of $S$ (*flats*) (see [Mo73]). For a given flat $f$ any basepoint $x \in f$ determines a decomposition of $f$ into *Weyl chambers* of $f$ based at $x$. Then $\mathcal{B}$ coincides with the space of asymptotic classes of Weyl chambers in $S$ (two chambers are asymptotic if they are within a bounded distance one from the other).

10.5. A flat with a distinguished class of asymptotic Weyl chambers is called an *oriented flat*. For an oriented flat $\overline{f}$ denote by $-\overline{f}$ the same flat with the orientation opposite to that of $\overline{f}$, and let $\pi_+(\overline{f}) \in \mathcal{B}$ (resp., $\pi_-(\overline{f}) = \pi_+(-\overline{f})$) be the corresponding asymptotic classes of Weyl chambers (the "endpoints" of $\overline{f}$). Denote by $\overline{f}_0$ the standard flat $f_0 = \exp(\mathfrak{A})o$ with the orientation determined by $\mathfrak{A}^+$. Let $b_0 = \pi_+(\overline{f}_0)$, and $b_w = \pi_+(w\overline{f}_0), w \in W$, where $W$ is the *Weyl group* which acts simply transitively on orientations of $f_0$. Denote by $w_0$ the element of $W$ (opposite to the identity) which is determined by the relation $w_0\overline{f}_0 = -\overline{f}_0$. Then the Bruhat decomposition of the group $\mathcal{G}$ and transitivity of the action of $\mathcal{G}$ on the space of oriented flats imply

THEOREM. *The $\mathcal{G}$-orbits $\mathcal{O}_w = \mathcal{G}(b_0, b_w), w \in W$ determine a stratification of the product $\mathcal{B} \times \mathcal{B}$, and $\mathcal{O}_{w_0}$ is the only orbit of maximal dimension. For any oriented flat $\overline{f}$ the pair of its endpoints $\bigl(\pi_-(\overline{f}), \pi_+(\overline{f})\bigr)$ belongs to $\mathcal{O}_{w_0}$, and conversely, for any pair $(b_-, b_+) \in \mathcal{O}_{w_0}$ there exists a unique oriented flat $\overline{f}(b_-, b_+)$ with endpoints $(b_-, b_+)$.*

*Remark.* In the rank 1 case flats are bilateral geodesics in $S$, and Weyl chambers are geodesic rays in $S$. The Weyl group consists of only 2 elements, and the orbits in $\mathcal{B} \times \mathcal{B}$ are the diagonal and its complement.

10.6. THEOREM ([GR85]). *Let $G$ be a discrete subgroup of a semi-simple Lie group $\mathcal{G}$, and $\mu$ be a probability measure on $G$ such that*

(i) *the semigroup $\mathrm{sgr}\,(\mu)$ generated by the support of $\mu$ contains a sequence $g_n$ such that $\langle \alpha, r(g_n o) \rangle \to \infty$ for any positive root $\alpha$,*

(ii) *no conjugate of the group $\mathrm{gr}\,(\mu)$ is contained in a finite union of left translations of degenerate double cosets from the Bruhat decomposition of $\mathcal{G}$.*

*Then*

(j) *for $\mathbf{P}$-a.e. sample path $\{x_n\}$ of the random walk $(G, \mu)$ the sequence $x_n o$ converges in the Satake-Furstenberg compactification of the symmetric space $S$,*

(jj) *the corresponding limit measure $\lambda$ is concentrated on the Furstenberg boundary $\mathcal{B}$, and it is the unique $\mu$-stationary measure on $\mathcal{B}$,* (jjj) *For any point $b_- \in \mathcal{B}$ the set $\{b_+ \in \mathcal{B} : (b_-, b_+) \in \mathcal{O}_{w_0}\}$ has full measure $\lambda$, where $\mathcal{O}_{w_0}$ is the maximal dimension stratum of the Bruhat stratification in $\mathcal{B} \times \mathcal{B}$ defined in* §10.5.

*Remark.* As it was noticed in [GM89], in the case when $\mathcal{G}$ is an algebraic group conditions (i) and (ii) follow from Zariski density of the semigroup sgr $(\mu)$ in $\mathcal{G}$. However, these conditions can be also satisfied without sgr $(\mu)$ being Zariski dense [GR89].

10.7. Conditions (i) and (ii) of Theorem 10.6 are clearly satisfied simultaneously for the measure $\mu$ and for the reflected measure $\check{\mu}$, and by Theorem 10.6 (jjj) the product $\lambda_- \times \lambda_+$ of the limit measures of the random walks $(G, \check{\mu})$ and $(G, \mu)$ is concentrated on the orbit $\mathcal{O}_{w_0}$. Since the flats in $S$ have polynomial growth, the strips in $G$ defined as

$$S(b_-, b_+) = \left\{g \in G : \operatorname{dist}(go, \overline{f}(b_-, b_+)) \leq R\right\},$$

where $\overline{f}(b_-, b_+)$ is the flat in $S$ with endpoints $b_-, b_+$, also have polynomial growth (and they are a.e. nonempty for a sufficiently large $R$). Theorem 6.5 (b) then implies

THEOREM. *Under conditions of Theorem* 10.6, *if the measure $\mu$ has finite first logarithmic moment $\sum \log \operatorname{dist}(go, o) \mu(g)$ and finite entropy $H(\mu)$, then the Poisson boundary of $(G, \mu)$ is nontrivial and isomorphic to the Furstenberg boundary $\mathcal{B}$ with the limit measure determined by Theorem* 10.6 (jj).

10.8. *Remarks.* 1. Theorem 10.3 was first announced in [Ka85]. For discrete subgroups of $\operatorname{SL}(d, \mathbb{R})$ another proof (under somewhat more restrictive conditions) was independently obtained in [Le85].

2. Conditions of Theorem 10.7 on the decay at infinity of the measure $\mu$ are more general than those of Theorem 10.3. As a trade-off, Theorem 10.7 requires irreducibility assumptions (i) and (ii) from Theorem 10.6, whereas Theorem 10.3 does not impose any conditions at all on the support of the measure $\mu$. Note that if the measure $\mu$ has a finite first moment, then under the conditions of Theorem 10.6 the vector $a(\mu)$ from Theorem 10.3 belongs to $\mathfrak{A}^+$ (see [GR85]), so that the orbit $\partial S_a$ is isomorphic to $\mathcal{B}$, and the descriptions of the Poisson boundary given in Theorems 10.3 and 10.7 coincide. Actually, Theorem 10.6 can also be used to identify the Poisson boundary for measures with a finite first moment without the irreducibility assumptions (i) and (ii). In this case instead of flats one has to take the symmetric subspaces of $S$ corresponding to pairs of boundary points which are not in general position with respect to the Bruhat decomposition and use the fact that the rate of escape along these subspaces is sublinear.

3. Our description of the Poisson boundary for Zariski dense subgroups (Theorem 10.7) coincides with the description of the Poisson boundary for absolutely continuous measures on semi-simple Lie groups obtained by Furstenberg [Fu63]. He proved that the Poisson boundary for an *arbitrary* initial distribution is a finite cover of $\mathcal{B}$ determined by periodicity properties of the measure $\mu$, and this cover is trivial for the initial distribution concentrated at the group identity (cf. §8.5).

4. The limit measure $\lambda$ on $\mathcal{B}$ does not have to be absolutely continuous with respect to the Haar measure on $\mathcal{B}$. Namely, for any finitely generated Zariski dense discrete subgroup $G \subset \mathcal{G}$ the author has constructed a symmetric finitely supported measure $\mu$ on $G$ with $\operatorname{gr}(\mu) = G$ such that $\lambda$ is singular (to be described elsewhere).

10.9. *Example.* Let $\mathcal{G} = \operatorname{SL}(d, \mathbb{R})$ with a maximal compact subgroup $\mathcal{K} = \operatorname{SO}(d)$. The map $g\mathcal{K} \mapsto \sqrt{gg^*}$ identifies the symmetric space $S = \mathcal{G}/\mathcal{K}$ with the set of positive definite $d \times d$ matrices with determinant 1, and the origin $o \cong \mathcal{K}$ corresponds to the identity matrix. Then the action of $\mathcal{G}$ on $S$ takes the form $(g, x) \mapsto \sqrt{gx^2g^*}$. Take for a Cartan subgroup $\mathcal{A} \subset \mathcal{G}$ the group of diagonal matrices with positive entries, so that its Lie algebra $\mathfrak{A}$ is the space $\{a = (\alpha_1, \alpha_2, \ldots, \alpha_d) \in \mathbb{R}^d : \sum \alpha_i = 0\}$, and choose a dominant Weyl chamber in $\mathfrak{A}$ as $\mathfrak{A}^+ = \{a \in \mathfrak{A} : \alpha_1 > \alpha_2 > \ldots \alpha_d\}$. The radial part $r(x) \in \overline{\mathfrak{A}}^+$ of a matrix $x \in S$ is the ordered vector of logarithms of its eigenvalues.

Geodesic rays in $S$ starting from $o$ have the form $\xi(t) = \xi_1^t$, where $\xi_1 \in S$ is a matrix at distance 1 from the origin $o$ (i.e., such that $\|r(\xi_1)\| = 1$), so that the visibility boundary $\partial S$ ($\equiv$ the space of geodesic rays issued from $o$) can be identified with the set $S_1$ of all such matrices, and a sequence $x_n \in S$ converges in the visibility compactification to $\xi_1 \in S_1 \cong \partial S$ if and only if $\log x_n / \|r(x_n)\| \to \log \xi_1$.

Matrices $\xi_1 \in S_1$ are parameterized by their eigenvalues and eigenspaces. However, it is more convenient to deal instead with the associated flags in $\mathbb{R}^d$. Namely, let $\lambda_1 > \cdots > \lambda_k$ be the distinct coordinates of the vector $r(\xi_1) = a$. Denote the eigenspace and the multiplicity of an eigenvalue $\lambda_i$ by $E_i \subset \mathbb{R}^d$ and $d_i = \dim E_i$, respectively. Then $\xi_1$ is uniquely determined by the vector $a$ and the flag $V_1 \subset V_2 \subset \cdots \subset V_k = \mathbb{R}^d$, where $V_i = \bigoplus_{j=k-i+1}^{k} E_j$. The spaces $V_i$ can be described by using the *Lyapunov exponents* $\chi(v) = \lim \log \|\xi_1^t v\|/t$, $v \in \mathbb{R}^d$ of the ray $\xi(t) = \xi_1^t$ as $V_i = \{v : \chi(v) \leq \lambda_{k-i+1}\}$ (here and below we assume $\chi(0) = -\infty$).

Thus, for a given vector $a \in \overline{\mathfrak{A}}_1^+$ the corresponding $\mathcal{G}$-orbit $\partial S_a \subset \partial S$ is the variety of flags in $\mathbb{R}^d$ of the type $(d_k, d_{k-1}+d_k, \ldots, d_2+d_3+\cdots+d_k)$, where $d_i$ are the multiplicities of the components of $a$. The Furstenberg boundary $\mathcal{B} = \mathcal{G}/\mathcal{P}$ of $S$ is isomorphic to nondegenerate orbits $\partial S_a$, $a \in \mathfrak{A}_1^+$ and coincides with the variety of full flags in $\mathbb{R}^d$, the minimal parabolic subgroup $\mathcal{P}$ being the group of upper triangular matrices.

For $\mathcal{G} = \mathrm{SL}(d, \mathbb{R})$ the first moment condition from Theorem 10.3 takes the form

$$\sum \log \|g\| \mu(g) < \infty , \tag{10.1}$$

and part (i) of the theorem is equivalent to saying that there exists a vector $a \in \overline{\mathfrak{A}}^+$ such that for **P**-a.e. sample path $\{x_n\}$ the sequence of matrices $x_n^*$ is *Lyapunov regular* with the *Lyapunov spectrum* $a$ (see [Ka89]). Namely, for any $v \in \mathbb{R}^d \setminus \{0\}$ there exists a limit $\chi(v) = \lim \log \|x_n^* v\|/n \in \{\lambda_1 > \cdots > \lambda_k\}$, and the subspaces $V_i = \{v \in \mathbb{R}^d : \chi(v) \leq \lambda_{k-i+1}\}$ have dimensions $\dim V_i = d_{k-i+1} + \cdots + d_k$, where $\lambda_i$ are the distinct components of $a$ with multiplicities $d_i$. If $a \neq 0$, then the limit of the sequence $\sqrt{x_n x_n^*}$ in the visibility compactification belongs to the orbit $\partial S_{a/\|a\|}$ and is determined by the *Lyapunov flag* $\{V_i\}$ of the sequence $x_n^*$. Therefore, Theorem 10.3 identifies the Poisson boundary for a measure $\mu$ on a discrete subgroup of $\mathrm{SL}(d, \mathbb{R})$ satisfying the moment condition (10.1) with the space of corresponding Lyapunov flags (the type of these flags is determined by the degeneracy of the Lyapunov spectrum).

The standard flat $f_0$ in $S$ is the set of diagonal matrices with positive entries, and the positive orientation on it determines the standard flag $b_0$ consisting of the subspaces $E_i \oplus \cdots \oplus E_d$, where $E_i$ are the coordinate subspaces of $\mathbb{R}^d$. The Weyl group $W$ is isomorphic to the symmetric group of the set $\{1, 2, \ldots, d\}$, and it acts on $f_0$ by permuting the diagonal entries. The element $w_0 \in W$ is the permutation $w_0 : (1, 2, \ldots, d-1, d) \mapsto (d, d-1, \ldots, 2, 1)$; the flag $b_{w_0} = w_0 b_0$ opposite to $b_0$ is obtained by reversing the order of coordinates and consists of subspaces $E_1 \oplus \cdots \oplus E_i$. For any vector $a \in \mathfrak{A}^+$ the matrices $\exp(ta) \in S$ converge in the Satake-Furstenberg compactification to $b_0$ (resp., to $b_{w_0}$) when $t \to \infty$ (resp., $t \to -\infty$). More generally, a pair of flags $(b_-, b_+)$ belongs to the $\mathcal{G}$-orbit of maximal dimension $\mathcal{O}_{w_0}$ in $\mathcal{B} \times \mathcal{B}$ if and only if there exists a matrix $g \in \mathcal{G}$ such that the sequence $g^n o = (g^n g^{*n})^{1/2}$ (resp., the sequence $g^{-n} o = (g^{-n} g^{*-n})^{1/2}$) converges in the Satake-Furstenberg compactification to $b_+$ (resp., $b_-$), i.e., if and only if the spectrum of $g$ is simple, absolute values of its eigenvalues are all pairwise distinct, and the Lyapunov flags of the sequences $g^{*n}$ and $g^{*-n}$ are $b_+$ and $b_-$, respectively. In fact, the stratification of $\mathcal{B}$ into the subvarieties $\{b_+ \in \mathcal{B} : (b_-, b_+) \in \mathcal{O}_w\}$ obtained for a fixed $b_- \in \mathcal{B}$ is the well known *Schubert stratification* of the flag variety.

Thus, Theorem 10.7 allows identification of the Poisson boundary with the flag variety for any measure $\mu$ on a discrete subgroup of $\mathrm{SL}(d, \mathbb{R})$ provided that $\mathrm{sgr}(\mu)$ is Zariski dense in $\mathrm{SL}(d, \mathbb{R})$, the measure $\mu$ satisfies the moment condition $\sum \log \log \|g\| \mu(g) < \infty$ and has a finite entropy $H(\mu)$.

CNRS UMR 6625, IRMAR, Campus Beaulieu, Rennes 35042, France
University of Manchester, Manchester M13 9PL, UK
*E-mail address*: Vadim.Kaimanovitch@univ-rennes1.fr